\documentclass{article}
\usepackage{mathptmx} 

\usepackage{latexsym,epsfig,bm}
\usepackage{amssymb}
\usepackage{color}
\usepackage{amsthm}
\usepackage{mathrsfs}

\DeclareSymbolFont{AMSb}{U}{msb}{m}{n}
\DeclareSymbolFontAlphabet{\mathbb}{AMSb}

\newcommand{\E}{{\cal E}}
\newcommand{\bS}{{\bf S}}
\newcommand{\bmGamma}{\bm\Gamma}
\newcommand{\bF}{{\bf F}}
\newcommand{\bA}{{\bf A}}
\newcommand{\toEF}{\stackrel{{\E}\sb F}\longrightarrow}

\def\xx{{j}}
\def\yy{{k}}

\newcommand\supp{\mathop{\rm supp}}
\newcommand\p{\partial}
\newcommand{\at}[1]{\vert\sb{\sb{#1}}}
\newcommand{\Spec}{\mathop{\rm Spec}}
\newcommand{\Bar}[1]{\mkern2mu\overline{\mkern-2mu#1\mkern-5mu}\mkern5mu}
\def\Re{{\rm Re\, }}
\def\Im{{\rm Im\,}}
\providecommand\C{{\mathbb C}}
\renewcommand\C{{\mathbb C}}
\newcommand{\R}{{\mathbb R}}

\newcommand{\Abs}[1]{\left\vert#1\right\vert}
\newcommand{\abs}[1]{\vert #1 \vert}
\newcommand{\Norm}[1]{\left\Vert #1 \right\Vert}
\newcommand{\norm}[1]{\Vert #1 \Vert}
\newcommand{\const}{{\rm const}}
\newcommand\cnst{\mathop{C}}
\newcommand\sothat{{\rm :}\ }
\newcommand\sgn{\mathop{\rm sgn}}

\newcommand\mod{\mathop{\rm mod}}

\providecommand{\ltor}[1]{
\ifnum #1=1{\it i}\else\ifnum #1=2{\it ii}\else\ifnum #1=3{\it iii}
\else\ifnum #1=4 {\it iv}\fi\fi\fi\fi
}

\DeclareMathSymbol{\varOmega}{\mathord}{letters}{"0A}

\font\thf cmssdc10 at 11pt
\theoremstyle{plain}
\newtheorem{theorem}{\thf Theorem}[section]

\newtheorem{lemma}[theorem]{\thf Lemma}
\newtheorem{corollary}[theorem]{\thf Corollary}
\newtheorem{proposition}[theorem]{\thf Proposition}

\theoremstyle{definition}
\newtheorem{definition}[theorem]{Definition}
\newtheorem{assumption}{Assumption}

\theoremstyle{remark}
\newtheorem{remark}[theorem]{Remark}

\makeatletter\@addtoreset{equation}{section}
\makeatother

\begin{document}
\title{
%
Global attractor for a nonlinear oscillator
coupled to the Klein-Gordon field
}

\author{
{\sc Alexander Komech}
\footnote{
On leave from Department of Mechanics and Mathematics,
Moscow State University, Moscow 119899, Russia.
Supported in part
by Max-Planck Institute for Mathematics in the Sciences (Leipzig),
the
Wolfgang Pauli Institute and the Faculty of Mathematics,
Vienna University,
by DFG Grant 436\,RUS\,113/615/0-1,
and by FWF Grant P19138-N13.}
\\
{\it Faculty of Mathematics, Wien A-1090, Austria}
\\ \\
{\sc Andrew Komech}
\footnote{
Supported in part
by Max-Planck Institute for Mathematics in the Sciences (Leipzig) and 
by the NSF Grant DMS-0434698.
}
\\
{\it 
Mathematics Department, Texas A\&M University,
College Station, TX, USA
}}

\maketitle
\begin{abstract}
The long-time asymptotics is analyzed for all finite energy
solutions to a model $\mathbf{U}(1)$-invariant
nonlinear Klein-Gordon equation in one dimension,
with the nonlinearity concentrated at a single point:
{\it each finite energy solution} converges
as $t\to\pm\infty$
to the set of all ``nonlinear eigenfunctions''
of the form $\psi(x)e\sp{-i\omega t}$.
The {\it global attraction}
is caused by the nonlinear
energy transfer from lower harmonics to the continuous spectrum
and subsequent dispersive radiation.

We justify this mechanism by the following novel
strategy
based on
\emph{inflation of spectrum
by the nonlinearity}.
We show that
any {\it omega-limit trajectory}
has the time-spectrum in the spectral gap $[-m,m]$ and
satisfies the original equation.
This equation implies the key
{\it spectral inclusion} for spectrum of the nonlinear term.
Then the application of the Titchmarsh Convolution Theorem
reduces
the spectrum of each omega-limit trajectory
to a single harmonic
$\omega\in[-m,m]$.

The research
is inspired by Bohr's postulate on
quantum transitions and
Schr\"o\-din\-ger's identification
of the quantum stationary states to the nonlinear eigenfunctions
of the coupled $\mathbf{U}(1)$-invariant
Maxwell-Schr\"o\-din\-ger and Maxwell-Dirac equations.
\end{abstract}


\section{Introduction}
\label{sect-introduction}
The long time asymptotics for
nonlinear
wave equations
have been the subject of intensive research,
starting with the pioneering papers by
Segal \cite{MR0153967,MR0152908},
Strauss \cite{MR0233062},
and Morawetz and Strauss \cite{MR0303097},
where the
nonlinear scattering and the local attraction to zero were considered.
The asymptotic stability of solitary waves
has been studied since the 1990s by
Soffer and Weinstein \cite{MR1071238,MR1170476},
Buslaev and Perelman \cite{MR1199635e,MR1334139},
and then by others.
The existing results
suggest that the set of
orbitally stable solitary waves
typically forms a local attractor,
that is, attracts
finite energy solutions
that were initially close to it.

\begin{figure}
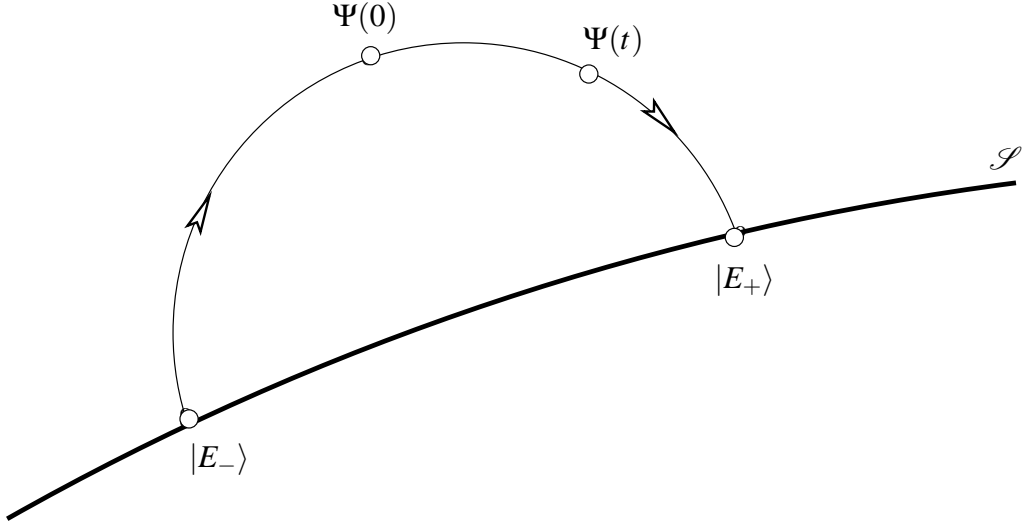

{
\input f2.pstex_t
}
\vskip 0.5cm
\caption{
Attraction of any trajectory $\Psi(t)$ to the set
of solitary waves
as $t\to\pm\infty$.}
\label{fig-gaa}
\end{figure}

In this paper, we consider the global attractor
for
all finite energy solutions.
For the first time, we prove that
in a particular $\mathbf{U}(1)$-invariant
dispersive Hamiltonian system
the global attractor is finite-dimensional
and is formed by solitary waves.
The investigation is inspired by Bohr's quantum transitions
(``quantum jumps").
Namely, according to Bohr's postulates \cite{bohr1913},
an unperturbed electron lives forever
in a \emph{quantum stationary state} $\vert E\rangle$
that has a definite value $E$ of the energy.
Under an external perturbation,
the electron
can
jump from one state to another:
\begin{equation}\label{transitions}
\vert E\sb{-}\rangle
\longmapsto
\vert E\sb{+}\rangle.
\end{equation}
The postulate
suggests the dynamical interpretation
of the transitions
as long-time attraction
\begin{equation}\label{ga}
\Psi(t)\longrightarrow\vert E\sb\pm\rangle,
\qquad
t\to\pm\infty
\end{equation}
for any trajectory $\Psi(t)$ of the corresponding dynamical system,
where the limiting states $\vert E\sb\pm\rangle$ generally
depend on the trajectory.
Then the quantum stationary states
should be viewed as the points of the \emph{global attractor}
$\mathcal{S}$ which is the set of all
limiting states
(see Figure~\ref{fig-gaa}).
Following de Broglie's ideas,
Schr\"o\-din\-ger identified the stationary states
$\vert E\rangle$
as the solutions
of the wave equation that have the form
\begin{equation}\label{solitary-waves-0}
\psi(x,t)=\phi\sb\omega(x)e\sp{-i\omega t},
\qquad
\omega=E/\hbar,
\end{equation}
where $\hbar$ is Planck's constant.
Then the attraction (\ref{ga})
takes the form of the long-time asymptotics
\begin{equation}\label{asymptotics}
\psi(x,t)
\sim
\psi\sb\pm(x,t)=\phi\sb{\omega\sb\pm}(x)e\sp{-i\omega\sb{\pm}t},
\qquad
t\to\pm\infty,
\end{equation}
that hold for each finite energy solution.
However, because of the superposition
principle,
the asymptotics of type (\ref{asymptotics})
are generally impossible for
the linear autonomous Schr\"o\-din\-ger equation
of type
\begin{equation}
(i\p\sb t-V(x))\psi(x,t)
=(-i\nabla-{\bm A}(x))\sp 2\psi(x,t),
\end{equation}
where $V(x)$ and ${\bm A}(x)$ are
scalar and vector potentials of a static external Maxwell field.
An adequate description of this process
requires to consider the Schr\"o\-din\-ger
(or Dirac)
equation
coupled to the Maxwell system
which governs the time evolution of the Maxwell 4-potential
$A(x,t)=(V(x,t),{\bm A}(x,t))$.
This coupling is inevitable indeed,
because, again by Bohr's postulates,
the transitions (\ref{transitions})
are followed by electromagnetic radiation responsible for the
atomic spectra.
The coupled Maxwell-Schr\"o\-din\-ger system
was initially introduced in \cite{Sch81109}.
It is a $\mathbf{U}(1)$-invariant
nonlinear Hamiltonian system.
Its global well-posedness was considered in \cite{MR1331696}.
One might expect the following
generalization of asymptotics (\ref{asymptotics})
for solutions to the coupled Maxwell-Schr\"o\-din\-ger
(or Maxwell-Dirac)
equations:
\begin{equation}\label{asymptotics-a}
(\psi(x,t), A(x,t))
\sim
\left(
\phi\sb{\omega\sb\pm}(x)e\sp{-i\omega\sb\pm t},
A\sb{\omega\sb\pm}(x)
\right),
\qquad
t\to\pm\infty.
\end{equation}
The asymptotics of this form are not available yet
in the context of coupled systems.
Let us mention that
the existence of the solitary  waves for
the coupled Maxwell-Dirac equations
was established in \cite{MR1386737}.

The asymptotics (\ref{asymptotics-a}) would mean that
the set of all solitary waves
\[
\{\left(\phi\sb\omega(x),A\sb\omega(x)\right):\omega\in\C\}
\]
forms a global attractor for the coupled system.
Similar convergence to a global attractor
is well known for dissipative systems, like Navier-Stokes equations
(see \cite{MR1156492,He81,MR1441312}).
In this context,
the global attractor is formed by the \emph{static stationary states},
and the corresponding asymptotics (\ref{asymptotics})
only hold for $t\to+\infty$ (and with $\omega\sb{+}=0$).

Our main impetus for writing this paper
was the natural question
whether dispersive Hamiltonian systems
could, in the same spirit,
possess finite dimensional global attractors,
and whether such attractors
are formed by the solitary waves.
We prove such a global attraction
for a model nonlinear Klein-Gordon equation
\begin{equation}\label{KG-0}
\ddot\psi(x,t)
=\psi''(x,t)-m^2\psi(x,t)+\delta(x)F(\psi(0,t)),
\qquad
\quad x\in\R.
\end{equation}
Here $m>0$, $\psi(x,t)$ is a continuous complex-valued wave function,
and $F$ is a nonlinearity.
The dots stand for the derivatives in $t$,
and the primes for the derivatives in $x$.
All derivatives and the equation are understood in
the sense of distributions.
Equation (\ref{KG-0}) describes the linear Klein-Gordon equation
coupled to the nonlinear oscillator.
We assume that equation (\ref{KG-0}) is $\mathbf{U}(1)$-invariant;
that is,
$$
F(e\sp{i\theta}\psi)=e\sp{i\theta}F(\psi),
\qquad
\psi\in\C,
\quad
\theta\in\R.
$$
Note that the group $\mathbf{U}(1)$ is also the (global) gauge group
of the coupled Maxwell-Schr\"o\-din\-ger and Maxwell-Dirac equations,
with the representation given by
\[
(\psi(x),A(x))\mapsto (e^{i\theta}\psi(x),A(x)).
\]
This gauge symmetry leads to the charge
conservation
and to the existence of the solitary wave solutions
of the form (\ref{asymptotics-a}) (see \cite{MR1386737}).
We clarify the special role of the
``nonlinear eigenfunctions'',
or solitary waves, of equation (\ref{KG-0})
which are finite energy solutions of type
(\ref{solitary-waves-0}):
\begin{equation}\label{soliton}
\psi\sb\omega(x,t)=\phi\sb\omega(x)e\sp{-i\omega t},\quad \omega\in\C.
\end{equation}
We prove that indeed they form the global attractor for all finite
energy solutions to (\ref{KG-0}).

Equation (\ref{KG-0}) has the following key features
of the
coupled Maxwell-Schr\"o\-din\-ger and Maxwell-Dirac
equations:
({\it i})
The linear part of this equation has a dispersive character;
({\it iii})
It is a nonlinear Hamiltonian system;
({\it iii})
It is $\mathbf{U}(1)$-invariant.
We suggest that just these features
are responsible for the global attraction,
such as
(\ref{asymptotics}), (\ref{asymptotics-a}),
to ``quantum stationary states''.
\bigskip

Let us introduce the set of all solitary waves.

\begin{definition} \label{dSS}
Let $\mathcal{S}$ be the set of all functions
$\phi\sb\omega(x)\in H\sp 1(\R)$ with $\omega\in\C$,
so that $\phi\sb\omega(x)e\sp{-i\omega t}$
is a solution to (\ref{KG-0}).
\end{definition}
Here $H\sp 1(\R)$ denotes the Sobolev space.
Generically,
the quotient $\mathcal{S}/\mathbf{U}(1)$
is isomorphic to a finite union of one-dimensional
intervals.
We will give an explicit construction
of the set of all solitary waves for equation (\ref{KG-0});
See Proposition~\ref{prop-solitons}
and its proof in Appendix~\ref{sect-solitons}.
Let us mention
that there are numerous results on the existence of
solitary wave solutions
of the form $\phi(x)e\sp{-i\omega t}$
to nonlinear Hamiltonian systems with $\mathbf{U}(1)$ symmetry
\cite{MR0454365,MR695535,MR695536,MR765961,MR847126,MR1344729}.
Typically, such solutions exist for
$\omega$ from an interval or a collection of intervals
of the real line.

\bigskip

Our main result is the following long-time asymptotics:
In the case when the nonlinearity $F$ is polynomial
of order strictly
greater than $1$,
we prove the attraction
of any finite energy solution
to the set $\mathcal{S}$ of all solitary waves:
\begin{equation}\label{attraction}
\psi(\cdot,t)\longrightarrow\mathcal{S},
\qquad t\to\pm\infty,
\end{equation}
where the convergence holds in local energy seminorms.
In the linear case, when $F(\psi)=a\psi$ with $a\in\R$,
there is generally no attraction
to $\mathcal{S}$; instead, we show that the global attractor
is
the linear span of all solitary waves, $ \langle\mathcal{S}\rangle$.
See Theorem~\ref{main-theorem-linear}.

\begin{remark}\label{remark-short}
Although we proved the attraction (\ref{attraction}) to $\mathcal{S}$,
we have not proved the attraction to a particular solitary wave,
falling short of proving (\ref{asymptotics}).
Hypothetically, a solution can be drifting along $\mathcal{S}$,
keeping asymptotically close to it,
but never approaching a particular solitary wave.
\end{remark}

\begin{remark}
The requirement that the nonlinearity $F$ is polynomial
allows us to apply the Titchmarsh convolution theorem
that is vital to the proof.
We do not know whether this requirement could be dropped.
\end{remark}

Let us mention related earlier results:

\begin{enumerate}
\item
The asymptotics of type (\ref{asymptotics})
were discovered first with $\psi\sb\pm=0$
in the scattering theory
\cite{MR0233062,MR0303097,MR498955,MR535231,MR654553,MR824083,MR1120284}.
In this case, the attractor $\mathcal{S}$ consists of the zero solution only,
and the asymptotics mean well-known \emph{local energy decay}.

\item
The \emph{global attraction}
of type (\ref{asymptotics})
with $\psi\sb\pm\ne 0$ and $\omega\sb{\pm}=0$
was established in
\cite{MR1203302e,MR1359949,MR1412428,MR1434147,MR1726676,MR1748357}
for a number of nonlinear
wave problems.
There the attractor $\mathcal{S}$
is the set of all \emph{static} stationary states.
Let us mention that this set could be infinite
and contain continuous components.

\item
First results on the asymptotics of type (\ref{asymptotics}),
with $\omega\sb\pm\ne 0$
were obtained for nonlinear $\mathbf{U}(1)$-invariant
Schr\"o\-din\-ger equations
in the context of asymptotic stability.
This establishes asymptotics of type (\ref{asymptotics})
but only for solutions close to the solitary waves,
proving the existence of a \emph{local attractor}.
This was first done in
\cite{MR1071238,MR1199635e,MR1170476,MR1334139},
and then developed in
\cite{MR1488355,MR1681113,MR1893394,MR1835384,MR1972870,MR2027616}
and other papers.
\end{enumerate}

The \emph{global attraction} (\ref{attraction})
to the solitary waves with  $\omega\ne 0$
was announced
for the first time
in \cite{MR2032730}  for equation (\ref{KG-0}).
In the present paper we give the detailed proofs,
and also add the well-posedness result which is not trivial
since the Dirac delta-function $\delta(x)$ does not belong to $L^2(\R)$.

Let us mention that the attraction (\ref{asymptotics})
for equation (\ref{KG-0}) with $m=0$
was proved in
\cite{MR1203302e,MR1359949}; In that case $\omega\sb\pm=0$.
Our proofs for $m>0$ are quite different from
\cite{MR1203302e,MR1359949}, and
are based
on a nonlinear spectral analysis
of \emph{omega-limit trajectories} for $t\to +\infty$
(and similarly for $t\to -\infty$).
First, we prove that their time-spectrum
is contained in a finite interval $[-m,m]$,
since the
spectral density
is absolutely continuous
for $\abs{\omega}>m$
and the corresponding component of the solution disperses completely.
Second,
the nonlinear equation (\ref{KG-0}) implies
the crucial spectral inclusion:
The nonlinearity does not inflate the spectrum
of any omega-limit trajectory.
Finally, the Titchmarsh convolution theorem
allows us to reduce
the spectrum of the omega-limit trajectory
to a single harmonic
$\omega\sb{+}\in[-m,m]$.
This implies the attraction (\ref{attraction}).

\begin{remark}
The global attraction
(\ref{asymptotics}), (\ref{asymptotics-a})
for $\mathbf{U}(1)$-invariant equations
suggests the corresponding extension
to general $\mathbf{G}$-invariant equations
($\mathbf{G}$ being the Lie group):
\begin{equation}\label{asymptotics-g}
\psi(x,t)
\sim
\psi\sb\pm(x,t)=e\sp{\bm\Omega\sb{\pm}t}\phi\sb\pm(x),
\qquad
t\to\pm\infty,
\end{equation}
where
$\bm\Omega\sb{\pm}$
belong to the corresponding Lie algebra
and
$e\sp{\bm\Omega\sb{\pm}t}$
are corresponding one-parameter subgroups.
Respectively,
the global attractor would consist of the solitary waves
(\ref{asymptotics-g}).
In particular, for the unitary group $\mathbf{G}=\mathbf{SU}(3)$,
the asymptotics (\ref{asymptotics-g})
relate the ``quantum stationary states'' to the structure of the
corresponding Lie algebra $\mathbf{su}(3)$.
On a seemingly related note,
let us mention that according to Gell-Mann -- Ne'eman theory \cite{GN64}
there is a correspondence
between the Lie algebras
and the classification of the elementary particles
which are the ``quantum stationary states''.
The correspondence
has been confirmed
experimentally by the discovery of the omega-minus Hyperon.
\end{remark}

\bigskip
The plan of the paper is as follows.
In Section~\ref{sect-results} we state the main assumptions and results.
Section~\ref{sect-splitting}
describes the exclusion of dispersive components
from the solution.
In Section~\ref{sect-spectral} we state the spectral properties of all
omega-limit trajectories
and apply the Titchmarsh Convolution Theorem.
For completeness,
we also give the exhaustive treatment of the linear case,
when
$F(\psi)=a\psi$ with $a\in\R$;
See Section~\ref{sect-linear-case}.
In Appendix~\ref{sect-solitons}, we collect the properties
of the solitary waves.
In Appendix~\ref{sect-quasimeasures}, we describe properties
of quasimeasures and corresponding multiplicators.
The global well-posedness of
equation (\ref{KG-0}) in $H\sp 1(\R)$
is proved in Appendix~\ref{sect-existence}.

\section{Main results}
\label{sect-results}

\subsection*{Model}
We consider the Cauchy problem for the Klein-Gordon equation
with the nonlinearity concentrated at a point:
\begin{equation}\label{KG}
\left\{
\begin{array}{l}
\ddot\psi(x,t)
=\psi''(x,t)-m^2\psi(x,t)+\delta(x)F(\psi(0,t)),
\qquad
x\in\R,
\quad t\in\R,
\\
\psi\at{t=0}=\psi\sb 0(x),
\qquad
\dot\psi\at{t=0}=\pi\sb 0(x).
\end{array}\right.
\end{equation}
If we identify a complex number $\psi=u+i v\in\C$
with the two-dimensional vector
$(u,v)\in\R\sp 2$,
then, physically, equation (\ref{KG}) describes small crosswise
oscillations of the infinite
string in three-dimensional space
$(x,u,v)$
stretched along the $x$-axis.
The string is subject to
the action of
an ``elastic force'' $-m^2\psi(x,t)$ and
coupled to a nonlinear oscillator
of force $F(\psi)$
attached at the point $x=0$.

We define
$\Psi(t)=
\left[\!\scriptsize{\begin{array}{c}
\psi(x,t)\\\pi(x,t)\end{array}}\!\right]$
and write the Cauchy problem
(\ref{KG})
in the vector form:
\begin{equation}\label{KG-cp}
\dot\Psi(t)
=
\left[\begin{array}{cc}0&1\\\p\sb x\sp 2-m^2&0\end{array}\right]
\Psi(t)
+
\delta(x)\left[\begin{array}{c}0\\F(\psi)\end{array}\right],
\qquad
\Psi\at{t=0}
=\Psi\sb 0
\equiv\left[\begin{array}{c}\psi\sb 0\\\pi\sb 0\end{array}\right].
\end{equation}
We will assume that
the oscillator force $F$ admits a real-valued potential:
\begin{equation}\label{P}
F(\psi)=-\nabla U(\psi),\quad\psi\in\C,
\qquad
U\in C\sp 2(\C),
\end{equation}
where the gradient is taken with respect to $\Re\psi$ and $\Im\psi$.
Then equation (\ref{KG-cp})
formally can be written as a Hamiltonian system,
\[
\dot\Psi(t)=J\,D\mathcal{H}(\Psi),
\qquad
J=\left[\begin{array}{cc}0&1\\-1&0\end{array}\right],
\]
where $D\mathcal{H}$ is the variational
derivative of the Hamilton functional
\begin{equation}\label{hamiltonian}
\mathcal{H}(\Psi)
=\frac 1 2
\int\limits\sb\R
\left(
\abs{\pi}\sp 2+\abs{\psi'}\sp 2+m^2\abs{\psi}\sp 2
\right)
dx
+U(\psi(0)),
\quad
\Psi=\left[\begin{array}{c}\psi(x)\\\pi(x)\end{array}\right].
\end{equation}
We
assume that the potential $U(\psi)$ is $\mathbf{U}(1)$-invariant,
where $\mathbf{U}(1)$ stands for the unitary group
$e\sp{i\theta}$, $\theta\in\R\mod 2\pi$:
Namely, we assume that
there exists $u\in C\sp 2(\R)$ such that
\begin{equation}\label{inv-u}
U(\psi)=u(\abs{\psi}\sp 2),
\qquad\psi\in\C.
\end{equation}

\begin{remark}
In the context of the model of the infinite string in $\R\sp{3}$
that we described after (\ref{KG}),
the potential $U(\psi)$
is rotation invariant with respect to the $x$-axis.
\end{remark}

Conditions (\ref{P}) and (\ref{inv-u})
imply that
\begin{equation}\label{def-a}
F(\psi)=\alpha(\abs{\psi}^2)\psi,
\qquad\psi\in\C,
\end{equation}
where
$\alpha(\cdot)=-2 u'(\cdot)\in C\sp 1(\R)$
is real-valued.
Therefore,
\begin{equation}\label{inv-f}
F(e\sp{i\theta}\psi)=e\sp{i\theta} F(\psi),
\qquad\theta\in\R,\quad\psi\in\C.
\end{equation}
Then the N\"other theorem formally
implies that
the functional
\begin{equation}\label{cal-Q}
\mathcal{Q}(\Psi)=\frac i2
\int\sb\R
\left(\overline\psi\pi-\overline\pi\psi\right)\,dx,
\qquad
\Psi=\left[\begin{array}{c}\psi(x)\\\pi(x)\end{array}\right],
\end{equation}
is
conserved for solutions $\Psi(t)$ to (\ref{KG-cp}).

We introduce
the phase space ${\E}$
of finite energy states for equation (\ref{KG-cp}).
Denote by $L\sp 2$ the complex Hilbert space $L\sp 2(\R)$
with the norm $\norm{\cdot}\sb{L\sp 2}$,
and denote by $\norm{\cdot}\sb{L\sp 2\sb R}$ the norm in $L\sp 2(-R,R)$
for $R>0$.

\begin{definition}
\begin{enumerate}
\item
${\E}$ is the Hilbert space of the states
$\Psi
=
\left[\begin{array}{c}\psi(x)
\\
\pi(x)
\end{array}\right]$,
with the norm
\begin{equation}\label{def-e}
\norm{\Psi}\sb{\E}^2
:=
\norm{ \pi}\sb{L\sp 2}^2
+\norm{\psi'}\sb{L\sp 2}^2+m^2\norm{\psi}\sb{L\sp 2}^2.
\end{equation}
\item
${\E}\sb F$ is the space ${\E}$ endowed
with the Fr\'echet topology defined by the seminorms
\begin{equation}\label{def-e-r}
\norm{\Psi }\sb{\E,R}^2
:=
\norm{ \pi}\sb{L\sp 2\sb R}^2
+
\norm{\psi'}\sb{L\sp 2\sb R}^2+m^2\norm{\psi}\sb{L\sp 2\sb R}^2,
\qquad
R>0.
\end{equation}
\end{enumerate}
\end{definition}

The
equation (\ref{KG-cp})
is formally a Hamiltonian system with
the phase space ${\E}$
and the Hamilton functional $\mathcal{H}$.
Both
$\mathcal{H}$
and $\mathcal{Q}$ are continuous functionals on ${\E}$.
Let us note that
${\E}={H\sp 1}\oplus L\sp 2$,
where $H\sp 1$ denotes the Sobolev space
\[
H\sp 1=H\sp 1(\R)
=\{\psi(x)\in L\sp 2(\R):\;\psi'(x)\in L\sp 2(\R)\}.
\]
We introduced into (\ref{def-e}), (\ref{def-e-r})
the factor $m^2>0$;
This provides the convenient relation
$\mathcal{H}(\Psi)=\frac 1 2\norm{\Psi}\sb{\E}^2+U(\psi(0))$.
The space  ${\E}\sb F$ is metrizable (but not complete).

\subsection*{Global well-posedness}

To have a priori estimates available for the proof of the global
well-posedness, we assume that
\begin{equation}\label{bound-below}
U(\psi)\ge {A}-{B}\abs{\psi}^2
\quad{\rm for}\ \psi\in\C,\quad
{\rm where}\ {A}\in\R\ {\rm and}\ 0\le {B}<m.
\end{equation}

\begin{theorem}\label{theorem-well-posedness}
Let $F(\psi)$ satisfy conditions (\ref{P}) and (\ref{inv-u}):
\[
F(\psi)=-\nabla U(\psi),\qquad
U(\psi)=u(\abs{\psi}^2),
\qquad u(\cdot)\in C\sp 2(\R).
\]
Additionally,
assume that (\ref{bound-below}) holds.
Then:
\begin{enumerate}
\item
For every $\Psi\sb 0\in {\E}$  the Cauchy problem
(\ref{KG-cp}) has a unique solution
$\Psi(t)\in C(\R,{\E})$.
\item
The map
$W(t):\;\Psi\sb 0\mapsto \Psi(t)$
is continuous in ${\E}$ and ${\E}\sb F$
for each $t\in\R$.
\item
The energy is conserved:
\begin{equation}\label{ec}
\mathcal{H}(\Psi(t))=\const,
\quad
t\in\R.
\end{equation}
\item
The following \emph{a priori} bound holds:
\begin{equation}\label{eb}
\norm{\Psi(t)}\sb{\E}
\le C(\Psi\sb 0),
\qquad t\in\R.
\end{equation}
\end{enumerate}
\end{theorem}

We prove this theorem in Appendix~\ref{sect-existence}.

\begin{remark}
The value of the charge is also conserved:
$\mathcal{Q}(\Psi(t))=\const$,
$\ t\in\R$.
\end{remark}

\subsection*{Solitary waves and the main theorem}

\begin{definition}\label{def-solitary-waves}
\begin{enumerate}
\item
The solitary waves of equation (\ref{KG-cp})
are solutions of the form
\begin{equation}\label{solitary-waves}
\Psi(t)=\Phi\sb\omega e\sp{-i\omega t},
\ {\rm where}\
\omega\in\C,
\ \Phi\sb\omega
=
\left[\!\begin{array}{c}\phi\sb\omega
\\
-i\omega\phi\sb\omega
\end{array}
\!\right],
\ \phi\sb\omega\in H\sp 1(\R).
\end{equation}
\item
The solitary manifold
is the set
$
\bS=
\left\{
\Phi\sb\omega\sothat\omega\in\C
\right\}
$
of all amplitudes $\Phi\sb\omega$.
\end{enumerate}
\end{definition}

Identity (\ref{inv-f}) implies that the set $\bS$
is invariant under multiplication by $e\sp{i\theta}$,
$\theta\in\R$.
Let us note that since $F(0)=0$ by (\ref{def-a}),
for any $\omega\in\C$
there
is a zero solitary wave with
$\phi\sb\omega(x)\equiv 0$.

Note that, according to (\ref{def-a}),
$\alpha(\abs{C}^2):=F(C)/C\in\R$
for any $C\in\C\backslash 0$.
We will need to distinguish the cases
when $F$ is linear and nonlinear; for this,
we introduce the following definition.

\begin{definition}\label{def-sn}
The function $F(\psi)$ is \emph{strictly nonlinear}
if the equation
$\alpha(C^2)=a$
has a discrete
(or empty)
set of positive roots
$C$ for each particular $a\in\R$.
\end{definition}

\begin{lemma}\label{lemma-omega-real}
If $F(\psi)$ is strictly nonlinear
in the sense of Definition~\ref{def-sn},
then
nonzero solitary waves
exist only for
$\omega\in\R$.
\end{lemma}

We prove this Lemma
in Appendix~\ref{sect-solitons}.

\begin{proposition}[Existence of solitary waves]
\label{prop-solitons}
Assume that $F(\psi)$ satisfies (\ref{inv-f}) and
that one of two following conditions holds:
\begin{enumerate}
\item
$F(\psi)$  is strictly nonlinear
in the sense of Definition~\ref{def-sn};
\item
$F(\psi)=a\psi$ with $a\in\R$.
\end{enumerate}
Then all nonzero solitary wave solutions
to {\rm (\ref{KG-cp})}
are given by {\rm (\ref{solitary-waves})} with
\begin{equation}\label{solitary-wave-profile}
\phi\sb\omega(x)=C e^{-\kappa\abs{x}},
\end{equation}
where $\kappa>0$, $\omega\in\C$, and $C\in\C\backslash 0$
satisfy the following relations:
\begin{equation}\label{kaka}
\alpha(\abs{C}^2)=2\kappa,
\qquad
\kappa^2=m^2-\omega^2.
\end{equation}
Additionally, if $F(\psi)$ is strictly nonlinear,
then $\omega\in(-m,m)$.
\end{proposition}
We prove this Proposition
in Appendix~\ref{sect-solitons}.

\begin{remark}
Let us denote
$\kappa\sb C={\alpha(\abs{C}^2)}/{2}$ and
$\omega\sb C\sp\pm=\pm\sqrt{m^2-\kappa\sb C\sp 2}$ for
$C \in\C$.
Then the relation  (\ref{kaka})
demonstrates that the set of all solitary waves
can be parametrized as follows:
\begin{enumerate}
\item
When $F(\psi)$ is strictly nonlinear,
in the sense of Definition~\ref{def-sn},
the profile function
$\phi\sb C(x)=C e^{i\theta}e^{-\kappa\sb C\abs{x}}$,
with $C\ge 0$ and
$\theta\in [0,2\pi]$,
corresponds to the solitary waves
with $\omega=\omega\sb C\sp\pm$
as long as $\kappa\sb C\in (0,m]$
(so that $\phi\sb C\in H\sp 1$
and $\omega$ is real in agreement with
Lemma~\ref{lemma-omega-real}).

\item
When $F(\psi)=a\psi$ with $a\in\R$,
we see from (\ref{kaka})
that $\kappa\sb C=a/2$
is constant.
If $a>0$,
the profile function
$\phi\sb C(x)=C e^{-a\abs{x}/2}$, with $C\in\C$,
corresponds to the solitary waves
with
$\omega=\pm\sqrt{m^2-\frac{a^2}{4}}$.
The restriction $\kappa\sb C\in (0,m]$
no longer applies
since the value of $\omega$ may be imaginary.
(This is different from the case of strictly nonlinear $F$,
when imaginary values of $\omega$
are prohibited by Lemma~\ref{lemma-omega-real}.)
If $a\le 0$,
then there is only the zero solitary wave solution.
\end{enumerate}

\end{remark}

As we mentioned before,
we need to assume that the nonlinearity is polynomial.
This assumption
is crucial in our argument:
It will allow to apply the Titchmarsh convolution theorem.
Now all our assumptions on $F$
can be summarized as follows.

\begin{assumption}\label{assumption-a}
\begin{equation}\label{f-is-such}
F(\psi)=-\nabla U(\psi),
\qquad
U(\psi)=\sum\limits\sb{n=0}\sp{N}u\sb n\abs{\psi}\sp{2n},
\end{equation}
where
$\ u\sb n\in\R$,
$\ u\sb N>0$,
$\ N\ge 2$.
\end{assumption}

This Assumption guarantees that the nonlinearity
$F$ satisfies
(\ref{P}) and (\ref{inv-u}),
and also the bound (\ref{bound-below})
from Theorem~\ref{theorem-well-posedness}.
Moreover, Assumption~\ref{assumption-a}
implies that $F$
is strictly nonlinear in the sense of Definition~\ref{def-sn}.
By Lemma~\ref{lemma-omega-real},
this in turn implies that
all nonzero solitary waves correspond to $\omega\in\R$.

Our main result is the following theorem.

\begin{theorem}[Main Theorem]
\label{main-theorem}
Let the nonlinearity $F(\psi)$
satisfy Assumption~\ref{assumption-a}.
Then for any $\Psi\sb 0\in \E$
the solution $\Psi(t)\in C(\R,{\E})$
to the Cauchy problem {\rm (\ref{KG-cp})}
with $\Psi(0)=\Psi\sb 0$
converges to $\bS$ in the space $\E\sb{F}$:
\begin{equation}\label{cal-A}
\Psi(t)\toEF \bS,
\quad t\to \pm\infty.
\end{equation}
\end{theorem}

Let us note that
the convergence to the set $\bS$
in the space $\E\sb{F}$ is equivalent to
\begin{equation}\label{metr}
\lim\sb{t\to\pm\infty}\rho(\Psi(t),\bS)=0,
\end{equation}
where $\rho$ is a
metric in the space $\E\sb{F}$
and
$\rho(\Psi(t),\bS):=\inf\limits\sb{\Phi\in \bS}\rho(\Psi(t),\Phi)$.

Let us also give the corresponding result for the linear case,
when $F(\psi)=a\psi$ with $a\in\R$.
We restrict our consideration to the case when
$a<2m$.
It is in this case that condition (\ref{bound-below}) is satisfied.
We do not consider the case
$a\ge 2m$, since in this case the solutions
are generally not bounded in $\E$-norm
(see Remark~\ref{remark-solitons-linear}),
while our arguments rely significantly on the bounds (\ref{eb}).
This case will be considered in more detail elsewhere.

\begin{theorem}[Linear case]
\label{main-theorem-linear}
Assume that $F(\psi)=a\psi$, where $a<2m$.
Then for any $\Psi\sb 0\in \E$
the solution $\Psi(t)\in C(\R,{\E})$
to the Cauchy problem {\rm (\ref{KG-cp})}
with $\Psi(0)=\Psi\sb 0$
converges in the space $\E\sb{F}$
to the linear span of $\bS$,
which we denote by $\langle\bS\rangle$:
\begin{equation}\label{cal-A-linear}
\Psi(t)\toEF\langle\bS\rangle,
\quad t\to \pm\infty.
\end{equation}

\end{theorem}

\begin{remark}

In Section~\ref{sect-linear-case}
we will show that:
\begin{enumerate}
\item
If $0<a<2m$,
then
$\langle\bS\rangle\ne\bS$.
Particular solutions
show that the attraction {\rm (\ref{cal-A})}
does not hold in general
(see Remark~\ref{noatt})
and has to be substituted by {\rm (\ref{cal-A-linear})}.
\item
If $a\le 0$,
$\langle\bS\rangle=\bS=\{0\}$.

\end{enumerate}
\end{remark}

\subsection*{Strategy of the proof}
For $m=0$
the global attraction of type {\rm (\ref{cal-A})}
is proved in \cite{MR1359949},
where the proof was based
on the direct calculation of the energy radiation
for the wave equation.
For the Klein-Gordon equation with $m>0$,
the dispersive relation $\omega\sp 2=k\sp 2+m^2$
results in the group velocities $v=\omega'(k)=k/\sqrt{k^2+m^2}$,
so every velocity $0\le\abs{v}<1$ is possible.
This complicates considerably the investigation of the energy propagation,
so
the approach \cite{MR1359949}
built on the fact that the group velocity was $\abs{v}=1$
no longer works.
To overcome this difficulty,
we introduce a new approach based on the
nonlinear spectral analysis
of the solution.

We prove the absolute continuity
of the spectrum of the solution  for $\abs{\omega}>m$.
This observation is similar to
the well-known Kato Theorem.
The proof is not obvious
and relies on the complex Fourier-Laplace transform and
the Wiener-Paley arguments.

We then split the solution into two components: Dispersive and bound,
with the frequencies  $\abs{\omega}>m$ and $\omega\in [-m,m]$,
respectively.
The dispersive component is an oscillatory integral of plane waves,
while the bound component
is a superposition of exponentially decaying functions.
The stationary phase argument leads to a local decay of
the dispersive component,
due to the absolute continuity of its spectrum.
This reduces the long-time behavior of the solution to the
behavior of the bound component.

Next, we establish the spectral representation
for the bound component. For this, we need to know
an optimal regularity of the corresponding spectral measure;
We have found out
that the spectral measure belongs to the space of
\emph{quasimeasures}
which are Fourier transforms of
bounded continuous functions,
\cite{Ga66}.
The spectral representation
implies compactness in the space of quasimeasures,
which in turn leads to the existence
of \emph{omega-limit trajectories}  for $t\to \infty$.

Further,  we prove that
an omega-limit trajectory itself
satisfies the nonlinear equation (\ref{KG-0}),
and this implies
the crucial spectral inclusion:
The spectrum of the nonlinear term
is included in the spectrum of the omega-limit trajectory.
We then reduce the spectrum of this limiting trajectory
to a single harmonic $\omega\sb{+}\in[-m,m]$
using the Titchmarsh convolution theorem
\cite{titchmarsh} (see also \cite[p.119]{MR1400006}
and \cite[Theorem 4.3.3]{MR1065136}).
In turn, this means that any omega-limit trajectory
lies in the manifold $\bS$ of the solitary waves,
which proves that $\bS$ is the global attractor.

Empirically,
the last part of our argument is a contemplation of
the radiative mechanism based on the
\emph{inflation of spectrum}
by the nonlinearity: A low-frequency perturbation
of the stationary state does not radiate
the energy until it generates (via a nonlinearity) ``a spectral line''
embedded in
the continuous spectrum outside $[-m,m]$.
This embedded spectral line
gives rise to the
wave packets
which bring the energy to infinity.
This radiative mechanism has been originally observed
in the numerical experiments with the nonlinear relativistic
Ginzburg-Landau equation (see \cite{MR2133795}).
The spectral inclusion for the  omega-limit trajectories
expresses their {\it nonradiative nature}:
The limiting trajectory cannot radiate
since the initial energy was bounded.

\section{Separation of dispersive components}
\label{sect-splitting}

It suffices to prove Theorem~\ref{main-theorem}
for $t\to+\infty$;
We will only consider the solution $\psi(x,t)$
restricted to $t\ge 0$.
In this section we
eliminate two
dispersive components from $\psi(x,t)$.

\subsection*{First dispersive component}
Let us split the solution
$\Psi(t)
=\left[\!\scriptsize{\begin{array}{c}\psi(x,t)\\\pi(x,t)\end{array}}\!\right]
$
into
$\Psi(t)=\Psi\sb 1(t)+\Psi\sb 2(t)$,
where
$
\Psi\sb 1(t)
=\left[\!\scriptsize{
\begin{array}{c}\psi\sb 1(x,t)\\\pi\sb 1(x,t)\end{array}}\!\right]
$
and
$
\Psi\sb 2(t)
=\left[\!\scriptsize{
\begin{array}{c}\psi\sb 2(x,t)\\\pi\sb 2(x,t)\end{array}}\!\right]
$
are defined for $t\ge 0$
as solutions to the following Cauchy problems:
\begin{eqnarray}
&&
\dot\Psi\sb 1(t)
=
\left[\begin{array}{cc}0&1\\\p\sb x\sp 2-m^2&0\end{array}\right]
\Psi\sb 1(t),
\qquad
\Psi\sb 1\at{t=0}
=\Psi\sb 0,
\label{KG-cp-1}
\\
\\
&&
\dot\Psi\sb 2(t)
=
\left[\begin{array}{cc}0&1\\\p\sb x\sp 2-m^2&0\end{array}\right]
\Psi\sb 2(t)
+
\delta(x)\left[\begin{array}{c}0\\f(t)\end{array}\right],
\qquad
\Psi\sb 2\at{t=0}=0,
\label{KG-cp-2}
\end{eqnarray}
where
$\Psi\sb 0
=\left[\!\scriptsize{\begin{array}{c}\psi\sb 0\\\pi\sb 0\end{array}}\!\right]$
is the initial data from
(\ref{KG-cp}),
and
\begin{equation}\label{def-f}
f(t):=F(\psi(0,t)),
\qquad
t\ge 0.
\end{equation}
Note that $\psi(0,\cdot)\in C\sb{b}(\Bar{\R\sp{+}})
$
by the Sobolev embedding
since $\Psi\in C\sb{b}(\Bar{\R\sp{+}}, \E)$
by Theorem~\ref{theorem-well-posedness}~({\it iv}).
Hence,  $f(\cdot)\in C\sb{b}(\Bar{\R\sp{+}})$.
On the other hand,
since $\Psi\sb 1(t)$
is a finite energy solution to the free Klein-Gordon equation,
we also have
\begin{equation}\label{psi-1-bounds}
\Psi\sb 1
\in C\sb{b}(\Bar{\R\sp{+}}, \E).
\end{equation}
Hence, the function $\Psi\sb 2(t)=\Psi(t)-\Psi\sb 1(t)$
also satisfies
\begin{equation}\label{psi-2-bounds}
\Psi\sb 2
\in C\sb{b}(\Bar{\R\sp{+}}, \E).
\end{equation}

\begin{lemma}\label{lemma-decay-psi1}
There is a local decay of $\Psi\sb 1$
 in ${\E}\sb F$ seminorms.
That is,
$\forall R>0$,
\begin{equation}\label{dp0}
\Norm{
\Psi\sb 1(t)
}\sb{{\E},R}
\to 0,
\qquad t\to\infty.
\end{equation}
\end{lemma}
\begin{proof}
We have to prove that
\begin{equation}
\norm{\Psi\sb 1(t)}\sb{\E,R}^2
=
\int\limits\sb{\abs{x}<R}\Big( \abs{\pi\sb 1(x,t)}\sp 2
+\abs{\psi\sb 1'(x,t)}\sp 2+m^2\abs{\psi\sb 1(x,t)}\sp 2\Big)\,dx
\end{equation}
goes to zero as $t$ tends to infinity.
Fix a cutoff function $\zeta(x)\in C\sb 0\sp\infty(\R)$
with $\zeta(x)=1$, $\abs{x}\le 1$ and $\zeta(x)=0$, $\abs{x}\ge 2$.
For $r>0$,
let $\zeta\sb r(x)=\zeta(x/r)$.
Denote by $\Phi\sb{r}(t)$ and $\Theta\sb{r}(t)$
the solutions to the free
Klein-Gordon equation with the initial data
$\zeta\sb r\Psi\sb 0$
and
$(1-\zeta\sb r)\Psi\sb 0$,
respectively,
so that $\Psi\sb 1(t)=\Phi\sb{r}(t)+\Theta\sb{r}(t)$.
Then there exists $C\sb r>0$ that depends on $r$
so that
$\norm{\Phi\sb{r}(t)}\sb{\E,R}^2
\le C\sb r (1+t)\sp{-1}$
for $t>0$,
since the solution
$\Phi\sb{r}$ is represented by the integral with the Green function,
which is the Bessel function decaying like
$(1+t)\sp{-1/2}$ (see e.g. \cite[$(2.7')$, Chapter I]{MR1364201}).
We then have:
\begin{equation}\label{ert}
\norm{\Psi\sb 1(t)}\sb{\E,R}^2
\le C\sb r(1+t)^{-1}
+
C
\norm{\Theta\sb{r}(t)}\sb{\E,R}^2
\end{equation}
where $r>0$ could be arbitrary.
To conclude that the left-hand side of (\ref{ert})
goes to zero,
it remains to note that
\begin{equation}\label{rec}
\norm{\Theta\sb{r}(t)}\sb{\E,R}
\le
\norm{\Theta\sb{r}(t)}\sb{\E}
=
\norm{\Theta\sb{r}(0)}\sb{\E}
\end{equation}
where the last relation
is
due to the energy conservation for the free Klein-Gordon equation,
and that the right-hand side of (\ref{rec})
could be made arbitrarily small
if $r>0$ is taken sufficiently large.
\qed\end{proof}

\subsection*{Complex Fourier-Laplace transform}

Let us analyze the complex Fourier-Laplace transform of
$\psi\sb 2(x,t)$:
\begin{equation}\label{FL}
\displaystyle
\tilde\psi\sb 2(x,\omega)
=\mathcal{F}\sp{+}\sb{t\to\omega}[\psi\sb 2(x,\cdot)]
:=
\int\sb 0\sp\infty e\sp{i\omega t}\psi\sb 2(x,t)\,dt,
\quad\omega\in\C\sp{+},
\end{equation}
where
$C\sp{+}:=\{z\in\C:\;\Im z>0\}$.
Due to (\ref{psi-2-bounds}),
$\tilde\psi\sb 2(\cdot,\omega)$
is an $H\sp 1$-valued analytic function of $\omega\in\C\sp{+}$.
Equation (\ref{KG-cp-2})
for $\psi\sb 2$ implies that
\[
-\omega\sp 2\tilde\psi\sb 2(x,\omega)
=
\tilde\psi\sb 2''(x,\omega)-m^2\tilde\psi\sb 2(x,\omega)
+\delta(x)\tilde f(\omega),\quad \omega\in\C\sp{+}.
\]
Hence,
the solution $\psi\sb 2(x,\omega)$ is a linear combination
of the fundamental solutions
which satisfy
\[
G\sb\pm''(x,\omega)+(\omega^2-m^2)G\sb\pm(x,\omega)
=\delta(x),
\qquad \omega\in\C\sp{+}.
\]
These solutions are given by
$\displaystyle
G\sb\pm(x,\omega)=\frac{e\sp{\pm ik(\omega)\abs{x}}}{\pm 2ik(\omega)}$,
where $k(\omega)$ stands for the analytic function
\begin{equation}\label{def-k}
k(\omega)=\sqrt{\omega\sp 2-m^2},
\qquad\Im k(\omega)>0,
\qquad
\omega\in\C\sp{+},
\end{equation}
which we extend
to $\omega\in\Bar{\C\sp{+}}$
by continuity.
We use the standard ``limiting absorption principle''
for the selection of the fundamental solution: Since
$\tilde\psi\sb 2(\cdot,\omega)\in H\sp 1$
for $\omega\in\C\sp{+}$,
only $G\sb{+}$ is appropriate,
because for $\omega\in\C\sp{+}$ the function $G\sb{+}(\cdot,\omega)$
is in $H\sp 1$ while $G\sb{-}$ is not.
Thus,
\begin{equation}\label{tilde-psi-tilde-f}
\tilde\psi\sb 2(x,\omega)=-\tilde f(\omega)G\sb{+}(x,\omega)
=-\tilde f(\omega)\frac{e\sp{i k(\omega)\abs{x}}}{2ik(\omega)},
\qquad\omega\in\C\sp{+}.
\end{equation}
Define
$\tilde z(\omega):=\mathcal{F}\sb{t\to\omega}\sp{+}[z(t)]$,
with
$z(t):=\psi\sb 2(0,t)$.
Then
$\tilde z(\omega)=-\tilde f(\omega)/(2ik(\omega))$, and
(\ref{tilde-psi-tilde-f}) becomes
\begin{equation}\label{p1}
\tilde\psi\sb 2(x,\omega)=\tilde z(\omega)e\sp{ik(\omega)\abs{x}},
\qquad\omega\in\C\sp{+}.
\end{equation}
Let us extend $\psi\sb 2(x,t)$ and $f(t)$ by zero for $t<0$:
\begin{equation}\label{zn}
\psi\sb 2(x,t)=0
\quad
{\rm and}
\quad
f(t)=0
\quad
{\rm for}\quad t<0.
\end{equation}
Then
\begin{equation}\label{psi2}
\psi\sb 2\in C\sb{b}(\R, H^1)
\end{equation}
by (\ref{psi-2-bounds}) since $\psi\sb 2(x,0+)=0$ by
initial conditions in (\ref{KG-cp-2}).
The Fourier transform
$\hat\psi\sb 2(\cdot,\omega):=\mathcal{F}\sb{t\to\omega}[\psi\sb 2(\cdot,t)]$
is a tempered $H\sp 1$-valued distribution of $\omega\in\R$
by (\ref{psi-2-bounds}).
The distribution $\hat\psi\sb 2(\cdot,\omega)$
is the boundary value of
the
analytic function $\tilde \psi\sb 2(\cdot,\omega)$,
in the following sense:
\begin{equation}\label{bvp1}
\hat\psi\sb 2(\cdot,\omega)
=\lim\limits\sb{\varepsilon\to 0+}\tilde\psi\sb 2(\cdot,\omega+i\varepsilon),
\qquad\omega\in\R,
\end{equation}
where the convergence
is in the space of
tempered distributions $\mathscr{S}'(\R\sb\omega,H\sp 1)$.
Indeed,
$\tilde\psi\sb 2(\cdot,\omega+i\varepsilon)
=\mathcal{F}\sb{t\to\omega}[\psi\sb 2(\cdot,t)e\sp{-\varepsilon t}]$
and
$\psi\sb 2(\cdot,t)e\sp{-\varepsilon t}
\mathop{\longrightarrow}\limits\sb{\varepsilon\to 0+}
\psi\sb 2(\cdot,t)$
where the convergence holds in $\mathscr{S}'(\R\sb t,H\sp 1)$ by (\ref{zn}).
Therefore, (\ref{bvp1}) holds by the continuity of the Fourier transform
$\mathcal{F}\sb{t\to\omega}$ in $\mathscr{S}'(\R)$.

Similarly to (\ref{bvp1}),
the distributions   $\hat z(\omega)$ and $\hat f(\omega)$,
$\omega\in\R$,
are the boundary values of the analytic in $\C\sp{+}$
functions $\tilde f(\omega)$ and $\tilde z(\omega)$,
$\omega\in\C\sp{+}$, respectively:
\begin{equation}\label{bv}
\hat z(\omega)=\lim\limits\sb{\varepsilon\to 0+}
\tilde z(\omega+i\varepsilon),
\qquad
\hat f(\omega)=\lim\limits\sb{\varepsilon\to 0+}
\tilde f(\omega+i\varepsilon), \quad \omega\in\R,
\end{equation}
since the functions
$z(t)$ and $f(t)$ are bounded for $t\ge 0$ and zeros for
$t<0$.
The convergence holds in the space of tempered distributions
$\mathscr{S}'(\R)$.

Let us justify that the representation (\ref{p1})
for $\hat\psi\sb 2(x,\omega)$
is also valid when $\omega\in\R$
if the multiplication in (\ref{p1})
is understood in the sense of quasimeasures
(see Appendix~\ref{sect-quasimeasures}).
\begin{proposition}\label{prop-uniform}
For any fixed $x\in\R$, the identity
\begin{equation}\label{p1r}
\hat\psi\sb 2(x,\omega)=
\hat z(\omega)e\sp{ik(\omega)\abs{x}},\quad\omega\in\R,
\end{equation}
holds in the sense of tempered distributions.
The
right-hand side is defined as the product of quasimeasure
$\hat z(\omega)$ by the multiplicator $e\sp{ik(\omega)\abs{x}}$.
\end{proposition}

\begin{proof}
The representation (\ref{p1r}) for $\omega\ne \pm m$
follows from (\ref{p1}) -- (\ref{bv}) since $e\sp{ik(\omega)\abs{x}}$
is a smooth function of $\omega\in\Bar{\C\sp{+}}$
outside the points $\pm m$.

So, we only need to justify
(\ref{p1}) in a neighborhood of each point $\omega=\pm m$.
The main problem is
the low regularity of $k(\omega)=\sqrt{\omega\sp 2-m^2}$
at the points $\pm m$.
Choose
a cutoff function $\zeta(\omega)\in C\sb 0\sp\infty(\R)$
such that
\begin{equation}\label{cfu}
\zeta\at{[-m-1,m+1]}\equiv 1.
\end{equation}
Let us note that $\psi\sb 2(x,\cdot)\in C\sb{b}(\R)$
for each particular $x\in\R$
by (\ref{psi2})
and the Sobolev embedding theorem.
Therefore,
for each $x\in\R$,
$\hat\psi\sb 2(x,\cdot)$
belongs to the space $\mathcal{QM}(\R)$
of quasimeasures
which are defined as functions with bounded continuous
Fourier transform
(see Definition~\ref{def-quasimeasure}).
In particular,
$\hat z(\cdot)$ is a quasimeasure since it is the Fourier transform
of the function $z(t):=\psi\sb 2(0,t)$.
On the other hand,
the function
$e\sp{i k(\omega)\abs{x}}\zeta(\omega)$ is a
multiplicator in $\mathcal{QM}(\R)$
by Lemma~\ref{lemma-quasi-2}~({\it i})
and  Lemma~\ref{lemma-quasi-1}~({\it i})
(see Appendix~\ref{sect-quasimeasures}).
Let us now prove that
\begin{equation}\label{p1rz}
\hat\psi\sb 2(x,\omega)\zeta(\omega)=
\hat z(\omega)e\sp{ik(\omega)\abs{x}}\zeta(\omega),\quad\omega\in\R
\end{equation}
in the sense of quasimeasures.
We define
$\mu\sb\varepsilon(\omega):=\tilde z(\omega+i\varepsilon)
=\mathcal{F}\sb{t\to\omega}[z(t)e\sp{-\varepsilon t}]$
for $\varepsilon>0$.
Then $\mu\sb\varepsilon(\omega)\in \mathcal{QM}(\R)$,
and $\mu\sb\varepsilon(\omega)\stackrel{\mathcal{QM}}
\longrightarrow \mu(\omega):= \hat z(\omega)$
as $\varepsilon\to 0+$
since
$z(t)e\sp{-\varepsilon t}\stackrel{C\sb{b,F}}\longrightarrow z(t)$
by (\ref{zn})
(see Definition~\ref{dA}).

Let us denote
$M\sb{x,\varepsilon}(\omega)
=
e\sp{i k(\omega+i\varepsilon)\abs{x}}\zeta(\omega)$
for $\omega\in\R$ and $\varepsilon\ge 0$.
By Lemmas~\ref{lemma-quasi-1}~({\it ii})
and~\ref{lemma-quasi-2}~({\it ii}),
$M\sb{x,\varepsilon}(\omega)$
are multiplicators in the space of quasimeasures.
This implies that
\begin{equation}\label{mve'}
\tilde z(\omega+i\varepsilon)
e\sp{i k(\omega+i\varepsilon)\abs{x}}
\zeta(\omega)
\stackrel{\mathcal{QM}}\longrightarrow \hat z(\omega)
e\sp{i k(\omega)\abs{x}}\zeta(\omega)
\quad{\rm as}\quad \varepsilon\to 0+.
\end{equation}
On the other hand, the left-hand side converges to the left-hand side
of (\ref{p1rz}) by
(\ref{p1}) and (\ref{bvp1}).
\qed\end{proof}

\subsection*{Absolutely continuous spectrum}
We study the regularity of the spectral density $\hat z(\omega)$
from (\ref{p1r}).
Denote
\begin{equation}\label{def-omega-delta}
\varOmega\sb\delta
:=(-\infty,-m-\delta)\cup(m+\delta,\infty),
\qquad
\delta\ge 0.
\end{equation}
Note that $\Bar{\varOmega\sb 0}=(-\infty,-m]\cup[m,\infty)$
coincides with the
continuous spectrum of the free Klein-Gordon equation,
and the function $\omega k(\omega)$ is positive for
$\omega\in\varOmega\sb{0}$.

\begin{proposition}\label{prop-continuity}
The distribution $\hat z(\omega)$ is absolutely continuous
for $\omega\in\varOmega\sb 0$,
and
$\hat z\in L\sp 1(\Bar{\varOmega\sb 0})$.
Moreover,
\begin{equation}\label{f2}
\int\sb{\varOmega\sb 0}\abs{\hat z(\omega)}\sp 2
\omega k(\omega)
\,d\omega<\infty.
\end{equation}
\end{proposition}
\begin{proof}
Let us first explain the main idea of the proof.
By (\ref{p1r}),
the function $\psi\sb{2}(x,t)$ \emph{formally} is a
{\it linear combination}
of the functions $e\sp{ik\abs{x}}$ with the
{\it amplitudes}
$\hat z(\omega)$:
\begin{equation}\label{osci}
\psi\sb{2}(x,t)=\frac 1{2\pi}\int\sb\R
\hat z(\omega)e\sp{ik(\omega)\abs{x}}e\sp{-i\omega t}\,d\omega,
\qquad x\in\R.
\end{equation}
For $\omega\in\varOmega\sb 0$,
the functions $e\sp{ik(\omega)\abs{x}}$ are of infinite $L^2$-norm,
while $\psi\sb{2}(\cdot,t)$ is of finite $L^2$-norm.
This is possible only if the amplitude is absolutely continuous
in $\varOmega\sb 0$:
For example, if we took $\hat z(\omega)=\delta(\omega-\omega\sb 0)$
with $\omega\sb 0\in\varOmega\sb 0$,
then $\psi\sb{2}(\cdot,t)$ would be of infinite $L^2$-norm.

For a rigorous proof, we use the Paley-Wiener arguments.
Namely, the Parseval identity and (\ref{psi-2-bounds})
imply that
\begin{equation}\label{PW}
\int\limits\sb\R
\norm{\tilde\psi\sb 2(\cdot,\omega+i\varepsilon)}\sb{H\sp 1}\sp 2\,d\omega
=2\pi\int\limits\sb 0\sp\infty e\sp{-2\varepsilon t}
\norm{\psi\sb 2(\cdot,t)}\sb{H\sp 1}\sp 2\,dt
\le\frac{\cnst}{\varepsilon},
\quad \varepsilon>0.
\end{equation}
On the other hand,
we can calculate the term in the left-hand side of (\ref{PW}) exactly.
First, according to (\ref{p1}),
\[
\tilde\psi\sb 2(\cdot,\omega+i\varepsilon)
=\tilde z(\omega+i\varepsilon)e^{ik(\omega+i\varepsilon)\abs{x}},
\]
hence (\ref{PW}) results in
\begin{equation}\label{PW-1}
\varepsilon
\int\sb\R
\abs{\tilde z(\omega+i\varepsilon)}^2
\norm{
e^{ik(\omega+i\varepsilon)\abs{x}}
}\sb{H\sp 1}\sp 2
\,d\omega
\le\cnst,
\qquad \varepsilon>0.
\end{equation}

Here is a crucial observation about the
norm of $e^{ik(\omega+i\varepsilon)\abs{x}}$.
\begin{lemma}
\begin{enumerate}
\item
For $\omega\in\R$,
\begin{equation}\label{wei}
\lim\sb{\varepsilon\to 0+}
\varepsilon \norm{e^{ik(\omega+i\varepsilon)\abs{x}}}\sb{H\sp 1}\sp 2
=n(\omega):=
\left\{
\begin{array}{ll}
\omega k(\omega),&\abs{\omega}>m\\
0,&\abs{\omega}<m
\end{array}
\right.,
\end{equation}
where the norm in $H\sp 1$ is chosen to be
$\norm{\psi}\sb{H\sp 1}
=\left(\norm{\psi'}\sb{L\sp 2}^2
+m^2\norm{\psi}\sb{L\sp 2}^2\right)\sp{1/2}.
$
\item
For any $\delta>0$
there exists $\varepsilon\sb\delta>0$ such that
for $\abs{\omega}>m+\delta$
and $\varepsilon\in(0,\varepsilon\sb\delta)$,
\begin{equation}\label{n-half}
\varepsilon \norm{e^{ik(\omega+i\varepsilon)\abs{x}}}\sb{H\sp 1}\sp 2
\ge n(\omega)/2.
\end{equation}
\end{enumerate}
\end{lemma}

\begin{proof}
Let us compute the $H\sp 1$-norm
using the Fourier space representation.
Setting $k\sb\varepsilon=k(\omega+i\varepsilon)$,
so that
$\Im k\sb\varepsilon>0$,
we get
$\mathcal{F}\sb{x\to k}
\left[e^{ik\sb \varepsilon\abs{x}}\right]
=2ik\sb\varepsilon/(k\sb \varepsilon^2-k^2)
$ for
$k\in\R$.
Hence,
\begin{equation}\label{g-h1-norm}
\norm{e^{ik\sb \varepsilon\abs{x}}}\sb{H\sp 1}\sp 2
=
\frac{2\abs{k\sb \varepsilon}^2}\pi
\int\sb{\R}
\frac{(k^2+m^2)dk}
{\abs{k\sb \varepsilon^2-k^2}^2}
=-4\Im\left[
\frac{(k\sb \varepsilon^2+m^2)\Bar{k\sb \varepsilon}}
{k\sb \varepsilon^2-\Bar{k\sb \varepsilon}^2}
\right].
\end{equation}
The last integral
is evaluated using the Cauchy theorem.
Substituting the expression $k\sb\varepsilon^2=(\omega+i\varepsilon)^2-m^2$,
we get:
\begin{equation}\label{norma}
\norm{e^{ik(\omega+i\varepsilon)\abs{x}}}\sb{H\sp 1}\sp 2
=\frac{1}{\varepsilon}
\Re
\left[
\frac{(\omega+i\varepsilon)^2\overline{k(\omega+i\varepsilon)}}{\omega}
\right],
\quad\varepsilon>0,
\ \omega\in\R,\ \omega\ne 0.
\end{equation}
The relation (\ref{wei}) follows
since the function $k(\omega)$ is real for $\abs{\omega}>m$,
but
is purely imaginary for  $\abs{\omega}<m$.

The second statement of the Lemma
follows
since $n(\omega)>0$ for $\abs{\omega}>m$, and
$n(\omega)\sim\abs{\omega}^2$ for $\abs{\omega}\to\infty$.
\qed\end{proof}

\begin{remark}
Note that $n(\omega)$ in (\ref{wei})
is zero for
$\abs{\omega}<m$,
since in that case the function
$e^{ik(\omega)\abs{x}}$ decays exponentially in $x$
and the $H\sp 1$-norm of
$e^{ik(\omega+i\varepsilon)\abs{x}}$
remains finite
when $\varepsilon\to 0+$.
\end{remark}

Substituting
(\ref{n-half})
into (\ref{PW-1}), we get:
\begin{equation}\label{fin}
\int\sb{\varOmega\sb\delta}
\abs{\tilde z(\omega+i\varepsilon)}\sp 2
\omega k(\omega)
\,d\omega
\le 2C,
\qquad
0<\varepsilon<\varepsilon\sb\delta,
\end{equation}
with the same $C$ as in (\ref{PW-1}),
and the region $\varOmega\sb\delta$ defined in
(\ref{def-omega-delta}).
We conclude that for each $\delta>0$ the set of functions
\[
g\sb{\delta,\varepsilon}(\omega)
=
\tilde z(\omega+i\varepsilon)
\abs{\omega k(\omega)}^{1/2},
\qquad
\varepsilon\in(0,\varepsilon\sb\delta),
\]
defined for $\omega\in\varOmega\sb\delta$,
is bounded in the Hilbert space $L\sp 2(\varOmega\sb\delta)$,
and, by the Banach Theorem, is weakly compact.
The convergence of the distributions (\ref{bv})
implies the following weak convergence in the Hilbert space
$L\sp 2(\varOmega\sb\delta)$:
\begin{equation}\label{wc}
g\sb{\delta,\varepsilon}
\rightharpoondown g\sb\delta,
\qquad \varepsilon\to 0+,
\end{equation}
where the limit function
$g\sb\delta(\omega)$ coincides with the distribution
$
\hat z(\omega)\abs{\omega k(\omega)}^{1/2}
$
restricted onto
${\varOmega\sb\delta}$.
It remains to note that the norms
of all
functions $g\sb\delta$, $\delta>0$,
are bounded in $L\sp 2(\varOmega\sb\delta)$ by (\ref{fin}),
hence (\ref{f2}) follows.
Finally, $\hat z(\omega)\in L\sp 1(\Bar{\varOmega\sb 0})$
by (\ref{f2}) and the Cauchy-Schwarz inequality.
\qed\end{proof}

Let us denote
\begin{equation}\label{dpo}
\hat z\sb{d}(\omega)
:=
\left\{
\begin{array}{ll}
\hat z(\omega),&\quad\omega\in\varOmega\sb 0,\\
0,&\quad\omega\in\R\setminus\varOmega\sb 0.
\end{array}\right.
\end{equation}
Then, by Proposition~\ref{prop-continuity},
$\hat z\sb{d}(\omega)\in L\sp 1(\R)$.

\subsection*{Second dispersive component}
Proposition~\ref{prop-continuity} and
the representation (\ref{osci})
suggest that we introduce the function
\begin{equation}\label{dd}
\psi\sb{d}(x,t)
=\frac 1{2\pi}\int\sb\R
\hat z\sb{d}(\omega)e\sp{ik(\omega)\abs{x}}e\sp{-i\omega t}\,d\omega,
\qquad x\in\R,\quad t\in\R,
\end{equation}
with $\hat z\sb{d}$ defined by (\ref{dpo}).
The Fourier transform of
$\psi\sb{d}(x,t)$
is given by the formula similar to (\ref{p1r}):
\begin{equation}
\label{ftd}
\hat\psi\sb{d}(x,\omega)
=
\displaystyle
\hat z\sb{d}(\omega)e\sp{ik(\omega)\abs{x}},\qquad x\in\R,
\qquad
\omega\in\R.
\end{equation}
We will show that $\psi\sb d(x,t)$
is a dispersive component of the solution $\psi(x,t)$,
in the following sense.
\begin{proposition}\label{prop-decay-psi-d}
$\psi\sb{d}(\cdot,t)$
is a bounded continuous $H\sp 1$-valued function:
\begin{equation}\label{bound-psi-d}
\psi\sb{d}(\cdot,t)\in C\sb{b}(\R,H\sp 1).
\end{equation}
The local energy decay holds
for $\psi\sb{d}(\cdot,t)$:
for any $R>0$,
\begin{equation}\label{decay-psi-d}
\Norm{
\left[\!{\scriptsize{
\begin{array}{c}
\psi\sb{d}(\cdot,t)\\\dot\psi\sb{d}(\cdot,t)
\end{array}
}}\!\right]
}\sb{{\E},R}
\to 0,
\qquad t\to\infty.
\end{equation}
\end{proposition}

\begin{proof}
Changing the variable,
we rewrite (\ref{dd}) as follows:
\begin{equation}\label{ppmr}
\psi\sb{d}(x,t)=\frac 1{2\pi}
\int\sb\R
\hat z(\omega(k))
e\sp{i k\abs{x}}e\sp{-i\omega(k)t}
\frac{k\,dk}{\omega(k)},
\qquad x\in\R,
\end{equation}
where $\omega(k)=\sqrt{k^2+m^2}$
is the branch analytic for $\Im k>0$ and continuous for  $\Im k\ge 0$.
Note that the function $\omega(k)$, $k\in\R\backslash 0$,
is the inverse function to
$k(\omega)$ defined on $\Bar{\C\sp{+}}$
(see (\ref{def-k}))
and restricted onto $\varOmega\sb 0=(-\infty,-m)\cup(m,\infty)$.
Let us introduce the functions
\begin{equation}\label{ppm}
\psi\sb\pm(x,t)\!:=\frac 1{2\pi}
\int\sb\R
\hat z(\omega(k))
e\sp{\pm i k x}e\sp{-i\omega(k) t}
\frac{kdk}{\omega(k)},
\qquad
x\in\R,
\quad
t\ge 0.
\end{equation}
Both functions $\psi\sb\pm(x,t)$
are solutions to the free Klein-Gordon equation
on the whole real line.
The (free Klein-Gordon) energy of each solution
is finite, since
\begin{eqnarray}
\Norm{
\left[\!{\scriptsize{
\begin{array}{c}\psi\sb\pm(\cdot,0)
\\\dot \psi\sb\pm(\cdot,0)\end{array}
}}\!\right]
}\sb{{\E}}^2
=
\int\sb{\R}
(\omega^2(k)+\abs{k}\sp 2+m^2)
\abs{\hat z(\omega(k))}^2
\frac{k^2}{\omega^2(k)}\,dk
\nonumber
\\
=
\int\sb{\R}
2\abs{\hat z(\omega(k))}^2
k^2\,dk
=
2\int\sb{\varOmega\sb 0}
\abs{\hat z(\omega)}^2
\omega k(\omega)\,d\omega
<\infty.
\nonumber
\end{eqnarray}
In the last inequality, we used (\ref{f2}).
Hence, both $\psi\sb{-}$ and $\psi\sb{+}$
are bounded continuous
$H\sp 1$-valued functions:
\begin{equation}\label{psipm}
\psi\sb\pm\in C\sb{b}(\R, H^1),
\end{equation}
and for any $R>0$
\begin{equation}\label{psi-pm-local-decay}
\Norm{
\left[\!{\scriptsize{
\begin{array}{c}\psi\sb\pm(\cdot,t)\\\dot \psi\sb\pm(\cdot,t)\end{array}
}}\!\right]
}\sb{{\E},R}
\to 0,
\qquad
t\to\infty
\end{equation}
by the same arguments as in
the proof of
Lemma~\ref{lemma-decay-psi1}.
The function
$\psi\sb{d}(x,t)$
coincides with $\psi\sb{+}(x,t)$ for $x\ge 0$
and with $\psi\sb{-}(x,t)$ for $x\le 0$:
\[
\psi\sb{d}(x,t)=\psi\sb{\pm}(x,t),
\qquad \pm x\ge 0.
\]
Moreover,
$\psi\sb{-}(0-,t)=\psi\sb{+}(0+,t)$,
so
$\psi\sb{d}(x,t)$ has no jump at $x=0$
and therefore
$\psi\sb{d}'(x,t)$ is square-integrable
over the whole $x$-axis.
Therefore,
 (\ref{bound-psi-d}) follows from (\ref{psipm}), and
 (\ref{decay-psi-d})
follows from (\ref{psi-pm-local-decay}).
\qed\end{proof}

\section{Bound component}
\label{sect-bound}

\subsection*{Spectral representation}
We introduce the bound component of the solution $\psi(x,t)$
by
\begin{equation}\label{bb}
\psi\sb{b}(x,t)
=\psi\sb{2}(x,t)-\psi\sb{d}(x,t)
=\psi(x,t)-\psi\sb 1(x,t)-\psi\sb{d}(x,t),
\ \ x\in\R,\ \ t\in\R.
\end{equation}
Then
(\ref{psi2}) and
Proposition~\ref{prop-decay-psi-d}
imply that
\begin{equation}\label{ebb}
\psi\sb{b}\in C\sb{b}(\R,H\sp 1).
\end{equation}
By (\ref{ebb}),
the function
\[
z\sb{b}(t):=\psi\sb{b}(0,t)
=\psi\sb 2(0,t)-\psi\sb{d}(0,t)
\]
is bounded and continuous.
Therefore, its Fourier transform
$
\hat z\sb{b}\in\mathscr{S}'(\R)
$
is a quasimeasure:
\begin{equation}\label{ftbs}
\hat z\sb{b}=\hat z-\hat z\sb{d}\in \mathcal{QM}(\R),
\qquad\supp \hat z\sb{b}\subset [-m,m],
\end{equation}
where the last inclusion follows from (\ref{dpo}).
Now (\ref{p1r}), (\ref{ftd}), and (\ref{bb})
imply the multiplicative relation
\begin{equation}\label{ftb0}
\hat\psi\sb{b}(x,\omega)
=
\hat z\sb{b}(\omega)e\sp{ik(\omega)\abs{x}}.
\end{equation}
We denote
\begin{equation}\label{def-kappa}
\kappa(\omega):=-i k(\omega)=\sqrt{m^2-\omega\sp 2},
\qquad
\Re\kappa(\omega)\ge 0
\quad{\rm for}\quad
\Im\omega\ge 0,
\end{equation}
where $k(\omega)$ was introduced in (\ref{def-k}).
Let us note that
$\Re\kappa(\omega)\ge 0$
and that $\kappa(\omega)>0$ for $\omega\in (-m,m)$.
We rewrite (\ref{ftb0}) as
\begin{equation}\label{ftb}
\hat\psi\sb{b}(x,\omega)
=
\hat z\sb{b}(\omega)e\sp{-\kappa(\omega)\abs{x}},
\qquad\omega\in\R.
\end{equation}
Therefore, (\ref{ftbs}) implies that
$\hat\psi\sb{b}(x,\omega)$ for any fixed $x\in\R$
is a quasimeasure with the support
$\supp\hat\psi\sb{b}(x,\cdot)\subset[-m,m]$,
and  finally,
\begin{equation}\label{bdef}
\psi\sb{b}(x,t)=\frac 1{2\pi}
\langle
\hat z\sb{b}(\omega)e\sp{-\kappa(\omega)\abs{x}},
e\sp{-i\omega t}\rangle,\qquad x\in\R,\quad t\in\R.
\end{equation}

\subsection*{Compactness}
We are going to prove a compactness of the set of
translations of the bound component,
$\{\psi\sb{b}(x,s+t)\sothat s\ge 0\}$.
We will derive the compactness from the following
uniform estimates which we deduce from
(\ref{ftbs}), (\ref{bdef})
by Lemma~\ref{lemma-quasi-1}.

\begin{proposition}
\begin{enumerate}
\item
The function $\psi\sb{b}(x,t)$ is smooth for $x\ne 0$
and $t\in\R$,
and
the following representation holds
for any fixed $\ x\ne 0$, $\ t\in\R$,
and any nonnegative integers $\xx$, $\yy$:
\begin{equation}\label{bqd}
\p\sb x\sp\xx\p\sb t\sp\yy
\psi\sb{b}(x,t)
=
\frac 1{2\pi}
\Big\langle \hat z\sb{b}(\omega)(-\kappa(\omega)\sgn x)^{\xx}
e\sp{-\kappa(\omega)\abs{x}},\,
(-i\omega)^{\yy} e\sp{-i\omega t}
\Big\rangle.
\end{equation}
\item
For any $R>0$,
there is a constant $C\sb{\xx,\yy,R}>0$ so that
\begin{equation}\label{bqda}
\sup\limits\sb{
0<\abs{x}\le R}
\,\,
\sup\limits\sb{t\in\R}
\abs{\p\sb x\sp\xx\p\sb t\sp\yy\psi\sb{b}(x,t)}
\le C\sb{\xx,\yy,R}.
\end{equation}
\end{enumerate}
\end{proposition}

\begin{remark}
Let us note that
the bounds (\ref{bqda}) are independent of $x$
and remain valid in the regions $x>0$ and $x<0$,
although the derivatives
$\p\sb x\sp\xx\p\sb t\sp\yy\psi\sb{b}(x,t)$
with $\xx\ne 0$
may have a jump at $x=0$.
(This is the case for the solitary waves in (\ref{solitary-wave-profile}).)
\end{remark}

\begin{proof}
({\it i})
\ The representation (\ref{bqd}) with $\xx=0$
and any $\yy$ follows
directly
from (\ref{ftbs}), (\ref{bdef}).
Further, consider, for example,
$\xx=1$ and $\yy=0$:
\begin{equation}\label{bqdx}
\p\sb x \psi\sb{b}(x,t)=\lim\limits\sb{\varepsilon\to 0}
\frac 1{2\pi}
\left\langle \hat z\sb{b}(\omega)
\frac{e\sp{-\kappa(\omega)\abs{x+\varepsilon}}-
e\sp{-\kappa(\omega)\abs{x}}}\varepsilon \zeta(\omega),
 e\sp{-i\omega t}\right\rangle
\end{equation}
if the limit exists.
Here $\zeta(\omega)\in C\sb{0}^\infty(\R)$ is any cutoff function
which satisfies (\ref{cfu}).
The relation (\ref{bqdx})
follows by Lemma~\ref{lemma-quasi-1}~({\it ii}),
if we verify that the following convergence
holds in $L\sp 1(\R\sb t)$:
\begin{equation}\label{cL1}
\mathcal{F}\sb{\omega\to t}^{-1}
\left[\frac{e\sp{-\kappa(\omega)\abs{x+\varepsilon}}
-e\sp{-\kappa(\omega)\abs{x}}}\varepsilon \zeta(\omega)\right]\stackrel{L\sp 1}\longrightarrow
\mathcal{F}\sb{\omega\to t}^{-1}
[\p\sb x e\sp{-\kappa(\omega)\abs{x}}\zeta(\omega)].
\end{equation}
We rewrite the expression in the brackets in the left-hand side
of (\ref{cL1}) as
\begin{equation}\label{cL1r}
\frac{e\sp{-\kappa(\omega)\abs{x+\varepsilon}}-
e\sp{-\kappa(\omega)\abs{x}}}\varepsilon
\zeta(\omega)
=
\int\sb 0\sp 1\p\sb x e\sp{-\kappa(\omega)\abs{x+\rho\varepsilon}}
\zeta(\omega)
\,d\rho.
\end{equation}
Now the convergence (\ref{cL1}) follows
from the Puiseux expansion of type (\ref{MHfA})
for $\p\sb x e\sp{-\kappa(\omega)\abs{x+\rho\varepsilon}}\zeta(\omega)$.

\noindent
({\it ii})
\ For fixed nonnegative integers $\xx$ and $\yy$, denote
\[
N\sb x(\omega)
=(-\kappa(\omega)\sgn x)\sp\xx
e\sp{-\kappa(\omega)\abs{x}}(-i\omega)\sp\yy\zeta(\omega).
\]
Lemma~\ref{lemma-quasi-2}~({\it iii})
implies that
$N\sb x(\omega)=\hat K\sb x(\omega)$,
with $K\sb x(\cdot)\in L\sp 1(\R)$.
Then (\ref{bqd}) becomes
\begin{equation}\label{bqdwb}
\p\sb x\sp\xx\p\sb t\sp\yy\psi\sb{b}(x,t)
=\frac 1{(2\pi)^2}
\int\sb\R z\sb{b}(t-\tau)K\sb x(\tau)d\tau.
\end{equation}
The bounds (\ref{bqda}) follow
by Lemma \ref{qmb} and (\ref{B*})
since
$z\sb{b}(t)=\psi\sb{b}(0,t)\in C\sb{b}(\R)$
by (\ref{ebb}).
\qed\end{proof}

\begin{corollary}\label{coco}
By the Ascoli-Arzel\`a Theorem,
for any sequence $s\sb{j}\to\infty$
there exists a subsequence $s\sb{j'}\to\infty$
such that
for any nonnegative integers $\xx$ and $\yy$,
\begin{equation}\label{olpd}
\p\sb x\sp\xx\p\sb t\sp\yy
\psi\sb{b}(x,s\sb{j'}+t)
\to
\p\sb x\sp\xx\p\sb t\sp\yy
\beta(x,t),
\qquad
x\ne 0,
\quad
t\in\R,
\end{equation}
for some
$\beta\in C\sb{b}(\R,H\sp 1)$.
The convergence in {\rm (\ref{olpd})}
is uniform in $x$ and $t$
as long as $\abs{x}+\abs{t}\le R$, for any $R>0$.
\end{corollary}

We call {\it omega-limit trajectory}
any function $\beta(x,t)$
that can appear as a limit in (\ref{olpd}).
Previous analysis demonstrates that the long-time asymptotics
of the solution $\psi(x,t)$ in $\E\sb F$
depends only on the bound component
$\psi\sb{b}(x,t)$.
By Corollary \ref{coco},
to conclude the proof of Theorem~\ref{main-theorem},
it suffices to check that every omega-limit trajectory
belongs to the set of solitary waves;
that is,
\begin{equation}\label{eidd}
\beta(x,t)
=
\phi\sb{\omega\sb{+}}(x)e\sp{-i\omega\sb{+}t},
\qquad
x,\,t\in\R,
\end{equation}
with some $\omega\sb{+}\in[-m,m]$.

\subsection*{Spectral identity for omega-limit trajectories}

Here we study the time spectrum of
the omega-limit trajectories.

\begin{definition}
Let $f$ be a tempered distribution.
By $\Spec f$
we denote the support of its Fourier transform:
\[
\Spec f:=\supp\tilde f.
\]
\end{definition}

\begin{proposition}\label{prop-beta}
\begin{enumerate}
\item
For any omega-limit trajectory
$\beta(x,t)$,
the following spectral representation
holds:
\begin{equation}\label{ber}
\beta(x,t)=\frac 1{2\pi}\langle \hat\gamma(\omega)e\sp{-\kappa(\omega)\abs{x}},e\sp{-i\omega t}\rangle,
\qquad x\in\R,\qquad t\in\R,
\end{equation}
where $\hat\gamma\in \mathcal{QM}(\R)$, and
\begin{equation}\label{spec-gamma}
\supp \hat\gamma\subset [-m,m].
\end{equation}
\item
The following bound holds:
\begin{equation}\label{ebbe}
\sup\limits\sb{t\in\R}
\norm{\beta(\cdot,t)}\sb{H\sp 1}<\infty.
\end{equation}
\end{enumerate}
\end{proposition}

Note that,
according to (\ref{ber}),
$\hat\gamma(\omega)$ is the Fourier
transform of the function $\gamma(t):= \beta(0,t)$, $t\in\R$.

\begin{proof}
The representation (\ref{bdef}) implies that
\begin{equation}\label{bqdwi}
\psi\sb{b}(x,s\sb{j}+t)
=
\frac{1}{2\pi}
\langle \hat z\sb{b}(\omega)
N\sb x(\omega)e\sp{-i\omega s\sb{j}},
e\sp{-i\omega t}\rangle,
\qquad x\ne 0,\qquad t\in\R,
\end{equation}
where $N\sb x$ corresponds to  $\xx=\yy=0$.
The convergence  (\ref{olpd}) and the bounds (\ref{bqda})
with  $\xx=\yy=0$ imply, by Definition~\ref{dA}, that
\begin{equation}\label{dztA}
z\sb{b}(s\sb{j'}+t)\stackrel{C\sb{b,F}}\longrightarrow \gamma(t),
\qquad s\sb{j'}\to\infty,
\end{equation}
where $\gamma(t)$, $t\in\R$,
is some continuous bounded function.
Hence, by Definition~\ref{dQ},
\begin{equation}\label{dztAd}
\hat z\sb{b}(\omega)e\sp{-i\omega s\sb{j'}}
\stackrel{\mathcal{QM}}\longrightarrow \hat\gamma(\omega),
\qquad
s\sb{j'}\to\infty.
\end{equation}
Now Lemma~\ref{lemma-quasi-1}~({\it ii})
and
Lemma~\ref{lemma-quasi-2}~({\it iii})
imply that
\begin{equation}\label{flA}
\hat z\sb{b}(\omega)
N\sb x(\omega)e\sp{-i\omega s\sb {k'}}
\stackrel{\mathcal{QM}}\longrightarrow \hat\gamma(\omega)N\sb x(\omega),
\qquad s\sb{j'}\to\infty.
\end{equation}
Hence, the representation (\ref{ber}) follows from
(\ref{bqdwi}), and (\ref{spec-gamma}) follows from (\ref{ftbs}).
Finally, the bound (\ref{ebbe})
follows from (\ref{ebb}) and (\ref{olpd}).
\qed\end{proof}

The relation (\ref{ber}) implies the
basic spectral identity:

\begin{corollary}
\label{corol}
For any omega-limit trajectory
$\beta(x,t)$,
\begin{equation}\label{SiI}
\Spec \beta(x,\cdot)=\Spec \gamma,
\qquad
x\in\R.
\end{equation}
\end{corollary}

\begin{remark}
It is mainly for the proof of (\ref{SiI})
that we develop the theory of the quasimeasures and multiplicators
in Appendix~\ref{sect-quasimeasures}.
This theory provides the compactness of the set of distributions
$\{
\hat z\sb{b}(\omega)e\sp{-i\omega s}\sothat s\ge 0\}$
in the space of  quasimeasures
(see (\ref{dztAd}))
and the spectral representation (\ref{ber}).
\end{remark}

\section{Nonlinear spectral analysis}
\label{sect-spectral}
Here we will derive (\ref{eidd})
from the following identity:
\begin{equation}\label{eide}
\gamma(t)=Ce^{-i\omega\sb{+}t},
\qquad
t\in\R,
\end{equation}
which will be proved in three steps.

\subsection*{Step 1}
The identity is equivalent to
$\hat\gamma(\omega)\sim\delta(\omega-\omega\sb{+})$, so we start with
an investigation of $\Spec \gamma:=\supp \hat \gamma$.
\begin{lemma}
For omega-limit trajectories
the following spectral inclusion holds:
\begin{equation}\label{Si}
\Spec F(\gamma(\cdot))\subset \Spec \gamma.
\end{equation}
\end{lemma}

\begin{proof}
The convergence  (\ref{olpd}) and equation (\ref{KG}),
together with Lemma~\ref{lemma-decay-psi1}
and Proposition~\ref{prop-decay-psi-d}~({\it ii}),
imply that the limiting trajectory $\beta(x,t)$ is a solution to
equation (\ref{KG}) (although $\psi\sb{b}(x,t)$ is not!):
\begin{equation}\label{KG-beta}
\ddot\beta(x,t)
=\beta''(x,t)-m^2\beta(x,t)+\delta(x)F(\beta(0,t)),
\qquad (x,t)\in\R\sp 2.
\end{equation}
Since $\beta(x,t)$ is smooth function for $x\le 0$ and $x\ge 0$,
we get
the following algebraic identity (cf. (\ref{ais})):
\begin{equation}
\label{AI}
0=\beta'(0+,t)-\beta'(0-,t)+F(\gamma(t)), \quad t\in\R.
\end{equation}
The identity implies the spectral inclusion
\begin{equation}\label{spectral-inclusion}
\Spec F(\gamma(\cdot))\subset
\Spec \beta'(0+,\cdot)\cup\Spec \beta'(0-,\cdot).
\end{equation}
On the other hand,
$\Spec \beta'(0+,\cdot)\cup\Spec \beta'(0-,\cdot)\subset\Spec \gamma$
by (\ref{SiI}).
Therefore, (\ref{spectral-inclusion})
implies
 (\ref{Si}).
\qed\end{proof}

\begin{remark}
The spectral inclusion (\ref{AI})
follows from the algebraic identity (\ref{AI}),
which in turn is a consequence of
the fact that $\beta(x,t)$ solves (\ref{KG}).
We cannot prove (\ref{spectral-inclusion})
for the function $\psi\sb{b}(x,t)$ since
generally it is \emph{not} a solution to (\ref{KG}).
\end{remark}

\subsection*{Step 2}
\begin{proposition}\label{pTN}
For any omega-limit trajectory,
the following identity holds:
\begin{equation}\label{C}
\abs{\gamma(t)}={\rm const},
\qquad t\in\R.
\end{equation}
\end{proposition}

\begin{proof}
We are going to show that (\ref{C})
follows from the key spectral relations (\ref{spec-gamma}), (\ref{Si}).
Our main assumption (\ref{f-is-such}) implies that the function
$F(t):=F(\gamma(t))$ admits the representation (cf. (\ref{def-a}))
\begin{equation}\label{frep}
F(t)=\alpha(t)\gamma(t),
\end{equation}
where, according to (\ref{f-is-such}),
\begin{equation}\label{arep}
\alpha(t)=-\sum\limits\sb{n=1}\sp{N}  2n u\sb n\abs{\gamma(t)}\sp{2n-2},
\qquad
N\ge 2;
\quad
u\sb N > 0.
\end{equation}
Both functions $\gamma(t)$ and $\alpha(t)$ are bounded
continuous functions in $\R$ by
Proposition~\ref{prop-beta}~({\it iii}).
Hence, $\gamma(t)$ and $\alpha(t)$ are
tempered distributions.
Furthermore, $\hat\gamma$ and $\hat{\overline\gamma}$ have the supports
contained in $[-m,m]$
by (\ref{spec-gamma}).
Hence, $\hat\alpha$ also has a bounded support
since it is
a sum of convolutions of finitely many
$\hat\gamma$ and $\hat{\overline\gamma}$ by (\ref{arep}).
Then the relation (\ref{frep})
translates into a convolution
in the Fourier space,
$
\hat F=\hat\alpha\ast\hat\gamma/(2\pi),
$
and the spectral inclusion (\ref{Si}) takes the following form:
\begin{equation}\label{Si-ast}
\supp\hat F=
\supp\,\hat\alpha\ast\hat\gamma\subset\supp\hat\gamma.
\end{equation}
Let us  denote $\bF=\supp \hat F$, $\bA=\supp \hat\alpha$, and
$\bmGamma=\supp\hat\gamma$.
Then the spectral inclusion (\ref{Si-ast}) reads as
\begin{equation}\label{Si-ast-r}
\bF\subset\bmGamma.
\end{equation}
On the other hand,
it is well known that
$\supp\hat\alpha\ast\hat\gamma\subset\supp\hat\alpha+\supp \hat\gamma$,
or
$
\bF \subset\bA+\bmGamma.
$
Moreover, the Titchmarsh convolution theorem
states that the last inclusion is exact for the ends
of the supports:

\begin{theorem}[The Titchmarsh Convolution Theorem]
{\it
Let $\hat\alpha$,
$\hat\gamma$ be two distributions in $\R$
with compact supports $\bA$ and $\bmGamma$ respectively,
and $\bF=\supp \hat\alpha*\hat\gamma$.
Then}
\begin{equation}\label{msc}
\inf\bF=\inf\bA
+\inf\bmGamma,
\qquad
\sup\bF=\sup\bA
+\sup\bmGamma.
\end{equation}
\end{theorem}

This theorem was proved first in \cite{titchmarsh}
for $\hat\alpha,\,\hat\gamma\in L\sp 1(\R)$ (see also \cite[p.119]{MR1400006}
and \cite[Theorem 4.3.3]{MR1065136}).
The Titchmarsh Convolution Theorem, together with (\ref{Si-ast-r}),
allows us to conclude that
$\inf\bA=\sup\bA=0$, hence
$
\bA\subset\{0\}.
$
Indeed,  (\ref{Si-ast-r}) and (\ref{msc})
result in
\begin{equation}\label{msci}
\inf\bF=\inf\bA
+\inf\bmGamma\ge\inf\bmGamma,
\qquad
\sup\bF=\sup\bA
+\sup\bmGamma\le \sup\bmGamma,
\end{equation}
so that $\inf\bA\ge 0\ge \sup\bA$.
Thus, we conclude that
$\supp\hat\alpha=\bA\subset\{0\}$,
therefore
the distribution
$\hat\alpha(\omega)$ is a finite
linear combination of $\delta(\omega)$ and its derivatives.
Then $\alpha(t)$ is a polynomial in $t$;
since $\alpha(t)$ is bounded by
Proposition~\ref{prop-beta}~({\it iii}),
we conclude that $\alpha(t)$ is constant.
Now the relation (\ref{C}) follows since
$\alpha(t)$ is a polynomial  in
$\abs{\gamma(t)}$, and its degree is strictly positive by
(\ref{arep}).
\qed\end{proof}
\begin{remark}
The boundedness of the spectrum of
both $\gamma(t)$
and $\alpha(t)$
is critical for our argument,
since otherwise the Titchmarsh convolution theorem
does not apply.
It is to ensure that the spectrum of
$\alpha(t)$ is also bounded
that we had to assume the polynomial
character of the nonlinearity
in Assumption~\ref{assumption-a}.
\end{remark}

\subsection*{Step 3}
Now the same Titchmarsh arguments imply that $\bmGamma:=\Spec \gamma$ is a point
$\omega\sb{+}\in[-m,m]$.
Indeed,
(\ref{C}) means that $\gamma(t) \overline\gamma(t)\equiv C$, hence
in the Fourier transform
$\hat\gamma \ast\hat{\overline\gamma}=2\pi C\delta(\omega)$.
Therefore,
if $\gamma$ is not identically zero,
the Titchmarsh Theorem implies that
\[
0=\sup\bmGamma+\sup(-\bmGamma)
=\sup\bmGamma-\inf\bmGamma.
\]
Hence
$\inf\bmGamma=\sup\bmGamma$
and therefore
$\bmGamma=\{\omega\sb{+}\}\subset [-m,m]$,
so that $\hat\gamma(\omega)$
is a finite linear combination
of $\delta(\omega-\omega\sb{+})$ and its derivatives.
As the matter of fact, the derivatives could not be present
because of the boundedness of
$\gamma(t)=\beta(0,t)$
that follows from
Proposition~\ref{prop-beta}~({\it iii}).
Thus,
$\hat\gamma\sim\delta(\omega-\omega\sb{+})$, which implies
(\ref{eide}).

\subsubsection*{Conclusion of the proof of Theorem~\ref{main-theorem}}
The representation
(\ref{ber}) implies that $\beta(x,t)
=C e\sp{-\kappa\sb{+}\abs{x}}e\sp{-i\omega\sb{+} t}$
since $\hat\gamma\sim\delta(\omega-\omega\sb{+})$.
Therefore, the equation
(\ref{KG-beta}) and the bound (\ref{ebbe})
imply
that $\beta(x,t)$
is a solitary wave.
This completes the proof of Theorem~\ref{main-theorem}.

\section{Linear case}
\label{sect-linear-case}
Let us now give a complete treatment of the linear case
and prove Theorem~\ref{main-theorem-linear}.
We assume that
$F(\psi)$ that enters (\ref{KG})
is given by
$F(\psi)=a\psi$,
where $a\in\R$, and $a<2m$.
Thus, the potential is given by
$U(\psi)=-a\abs{\psi}^2/2$,
and
we consider
the equation
\begin{equation}
\ddot\psi(x,t)
=\psi''(x,t)-m^2\psi(x,t)+a\delta(x)\psi(0,t),
\qquad
x\in\R,
\quad t\in\R.
\label{KG-linear}
\end{equation}
All  conclusions of
Theorem~\ref{theorem-well-posedness}
on global well-posedness
hold for equation (\ref{KG-linear})
with $a<2m$ since in this case
the condition (\ref{bound-below}) is satisfied.
Let us note that if $a\ge 2m$,
then the conclusions ({\it i}), ({\it ii}),
and ({\it iii}) of Theorem~\ref{theorem-well-posedness},
are still valid
(their proofs in  Appendix~\ref{sect-existence}
apply for bounded times,
and then the conclusions
follow
for all times by the linearity of the equation).
On the other hand,
the a priori bound (\ref{eb}) is generally violated
when $a/2\ge m$ (see below).

\begin{remark}\label{remark-solitons-linear}
Let us summarize the properties of the solitary
waves for the linear case that follow from Proposition~\ref{prop-solitons}.
Note that, according to (\ref{kaka}), $\kappa=a/2$.
\begin{enumerate}
\item
For $a\le 0$ there are no nonzero solitary waves
since we need $\kappa>0$ for (\ref{solitary-wave-profile})
to be from $H\sp 1$.
\item
When $a>0$, $a\ne 2m$,
all solitary waves are given by
$
\phi\sb\omega(x)=C e^{-a\abs{x}/2},
$
where $C\in\C$ and
$\omega=\pm\omega\sb a$,
where $\omega\sb a:=\sqrt{m^2-{a^2}/{4}}$.
Note that
if $a>2m$,
then the values $\pm\omega\sb a$ are purely imaginary
and the corresponding solitary waves are exponentially growing.
\item
If $a=2m$, then $\omega\sb 0=0$ and
there is a nonzero static solitary wave solution
$\phi\sb 0(x)=e^{-m\abs{x}}$.
Besides,
there is secular (linearly growing) solution
$te^{-m\abs{x}}$.

\end{enumerate}

\end{remark}

\begin{remark}
When $a>2m$,
the values of $\omega$ are purely imaginary,
and the $\E$-norm of solitary waves
that correspond to $\pm\Im\omega>0$
grows exponentially for $t\to\pm\infty$.
When $a=2m$, we have $\omega=0$;
The $\E$-norm
of the secular solution
grows linearly in time.
In both cases ($a\ge 2m$),
the a priori bound (\ref{eb}) fails.
This illustrates that
condition (\ref{bound-below}) is sharp,
since this condition fails for
the potential
$U(\psi)=-\displaystyle\frac{a}{2}\abs{\psi}^2$
with $a\ge 2m$.
\end{remark}

\subsubsection*{Proof of Theorem~\ref{main-theorem-linear}}

Let us prove the global attraction
to the set $ \langle\bS\rangle$.
We proceed as in the proof of Theorem~\ref{main-theorem}
until we get to equation (\ref{KG-beta}).
Since now $F(\psi)=a\psi$,
(\ref{KG-beta})
takes the following form:
\begin{equation}\label{KG-beta-linear}
\ddot\beta(x,t)
=\beta''(x,t)-m^2\beta(x,t)
+a\delta(x)\beta(0,t),
\qquad (x,t)\in\R\sp 2.
\end{equation}
Now we cannot use  the Titchmarsh arguments,
and we have to solve the equation  directly to prove
that
\begin{equation}\label{fol}
\left[\begin{array}{c}
\beta(\cdot,t)
\\
\dot\beta(\cdot,t)
\end{array}
\right]
\in \langle\bS\rangle
\qquad
{\rm for}\quad t\in\R.
\end{equation}
In
the Fourier transform
$\hat\beta(x,\omega)=\mathcal{F}\sb{t\to\omega}[\beta(x,t)]$
the equation (\ref{KG-beta-linear})  becomes
\begin{equation}\label{KG-beta-linear-f}
-\omega^2\hat\beta(x,\omega)
=\hat\beta''(x,\omega)-m^2\hat\beta(x,\omega)
+a\delta(x)\hat\beta(0,\omega),
\qquad (x,\omega)\in\R\sp 2.
\end{equation}
On the other hand, the representation
(\ref{ber}) implies that
\begin{equation}\label{beri}
\hat\beta(x,\omega)=\hat\gamma(\omega)e^{-\kappa(\omega)\abs{x}}
\end{equation}
Substituting this into  (\ref{KG-beta-linear-f}), we obtain
\begin{equation}\label{berig}
2\kappa(\omega)\hat\gamma(\omega)\delta(x)=a\delta(x)\hat\gamma(\omega).
\end{equation}
Therefore,
on the support of the distribution $\hat\gamma(\omega)$,
the identity holds
\begin{equation}\label{berigk}
2\kappa(\omega)=a,
\end{equation}
hence
$\supp \hat\gamma\subset \bmGamma\sb a:=\{\omega\in [-m,m]: 2\kappa(\omega)=a\}$
by (\ref{spec-gamma}).
Now let us consider two cases.

\begin{enumerate}
\item
In the case $0<a<2m$,
according to Remark~\ref{remark-solitons-linear}~({\it ii}),
the set of finite energy solitary waves is given by
\begin{equation}\label{ssw-1}
\bS=\left\{C\sb 1
\left[\begin{array}{c}e^{-a\abs{x}/2}\\
i\omega\sb a e^{-a\abs{x}/2}
\end{array}
\right]
+
C\sb 2\left[\begin{array}{c}e^{-a\abs{x}/2}\\
-i\omega\sb a e^{-a\abs{x}/2}
\end{array}\right]
\sothat\;
C\sb 1,\,\,C\sb 2\in\C
\right\}.
\end{equation}

On the other hand, the set
$\bmGamma\sb a$
contains exactly two points $\pm\omega\sb a$
since $0<a<2m$. Hence, $\hat\gamma$ is a linear combination of
$\delta(\omega\pm\omega\sb a)$ and their derivatives.
The derivatives are forbidden since
$\gamma(t)$ is bounded, so finally
\begin{equation}\label{fina}
\beta(x,t)
=
\left(
C\sb 1 e^{i\omega\sb a t}+C\sb 2 e^{-i\omega\sb a t}
\right)
e^{-a\abs{x}/2}.
\end{equation}
Now  (\ref{fol}) follows from  (\ref{fina}).
\item
In the case $a\le 0$,
the set of finite energy solitary waves consists of
the zero solution only by Remark~\ref{remark-solitons-linear}~({\it i}).
For $a<0$, the set $\bmGamma\sb a$ is empty,
hence $\beta(x,t)=0$ and (\ref{fol}) follows.
When $a=0$, we have
$\omega\sb a=m$
and
$\bmGamma\sb a=\{-m\}\cup\{m\}$.
Any omega-limit point
$\beta$ is given by
(\ref{fina})
with $a=0$.
Since
$\beta(\cdot,t)\in H\sp 1$,
we conclude that $C\sb 1=C\sb 2=0$ in (\ref{fina}),
so that $\beta(x,t)=0$ and the inclusion (\ref{fol}) follows.
\end{enumerate}

This finishes the proof of Theorem~\ref{main-theorem-linear}.

\begin{remark}\label{noatt}
For $0<a<2m$,
a particular exact solution to (\ref{KG-linear}),
e.g.  (\ref{fina}),
with $C\sb 1\ne 0$ and $C\sb 2\ne 0$,
shows that in general there could be no attraction to
$\bS$.
\end{remark}

\appendix

\section{Appendix: Solitary waves}
\label{sect-solitons}

Here we prove Lemma~\ref{lemma-omega-real}
and Proposition~\ref{prop-solitons}.

\subsubsection*{Proof of Lemma~\ref{lemma-omega-real}}
Substituting $\phi\sb\omega(x) e^{-i\omega t}$
into (\ref{KG}),
we get the equation
\begin{equation}\label{stat-eqn}
-\omega^2\phi\sb\omega(x)e^{-i\omega t}
=
\phi\sb\omega''(x)e^{-i\omega t}
-m^2\phi\sb\omega(x)e^{-i\omega t}
+\delta(x)F(e^{-i\omega t}\phi\sb\omega(0)),
\end{equation}
where $(x,t)\in\R\times\R$.
We can assume that $\phi\sb\omega(0)\ne 0$.
Indeed, if $\phi\sb\omega(0)=0$,
then (\ref{stat-eqn})
turns into a homogeneous second-order linear differential equation,
which together with the inclusion $\phi\sb\omega\in H\sp 1(\R)$
results in $\phi\sb\omega(x)\equiv 0$.
Equation (\ref{stat-eqn}) leads to the identity
$
e^{-i\omega t}\Delta
=
F(e^{-i\omega t}\phi\sb\omega(0)),
$
with
$\Delta=\phi\sb\omega'(0-)-\phi\sb\omega'(0+)$.
This results in
\begin{equation}\label{delta-f-phi}
\frac{\Delta}{\phi\sb\omega(0)}
=\frac{F(e^{-i\omega t}\phi\sb\omega(0))}{e^{-i\omega t}\phi\sb\omega(0)}
=\frac{F(e^{t\,\Im\omega}\phi\sb\omega(0))}{e^{t\,\Im\omega}\phi\sb\omega(0)},
\qquad
t\in\R.
\end{equation}
We used (\ref{inv-f}) in the last equality.
The condition that $F(\psi)$ is strictly nonlinear
(in the sense of Definition~\ref{def-sn})
implies that (\ref{delta-f-phi})
only holds at discrete values  of $t\,\Im\omega$;
thus, $\Im\omega=0$, finishing the proof.

\subsection*{Proof of Proposition~\ref{prop-solitons}}
The relation (\ref{stat-eqn}) turns into
the following eigenvalue problem:
\begin{equation}\label{NEP}
-\omega\sp 2\phi\sb\omega(x)
=
\phi\sb\omega''(x)-m^2\phi\sb\omega(x)
+\delta(x)F(\phi\sb\omega(x)),\quad x\in\R.
\end{equation}
The phase factor $e^{-i\omega t}$
can be canceled out
because either $F(\psi)=a\psi$ or,
when $F$ is strictly nonlinear,
we can use (\ref{inv-f})
(since in this case $\omega\in\R$ by Lemma~\ref{lemma-omega-real}).
Equation (\ref{NEP}) implies that away from the origin
we have
$$
\phi\sb\omega''(x)=(m^2-\omega^2)\phi\sb\omega(x),
\qquad
x\ne 0,
$$
hence
$\phi\sb\omega(x)=C\sb\pm e\sp{-\kappa\sb\pm\abs{x}}$ for $\pm x>0$,
where $\kappa\sb\pm$
satisfy
$\kappa\sb\pm^2=m^2-\omega\sp 2$.
Since $\phi\sb\omega(x)\in H\sp 1$,
it is imperative that $\kappa\sb\pm>0$;
we conclude that
$\abs{\omega}<m$
and that
$\kappa\sb\pm=\sqrt{m^2-\omega^2}>0$.
Moreover, since
the function $\phi\sb\omega(x)$ is continuous,
$C\sb{-}=C\sb{+}=C\ne 0$
(since we are looking for nonzero solitary waves).
We see that
\begin{equation}\label{profile}
\phi\sb\omega(x)=C e\sp{-\kappa\abs{x}},
\qquad
C\ne 0,
\qquad
\kappa\equiv\sqrt{m^2-\omega^2}>0.
\end{equation}
Equation (\ref{NEP}) implies the following gluing condition at $x=0$:
\begin{equation}\label{ais}
0=\phi\sb\omega'(0+)-\phi\sb\omega'(0-)+F(\phi\sb\omega(0)).
\end{equation}
This condition and (\ref{profile})
lead to the equation
$
2\kappa C=F(C)
$
which is equivalent to (\ref{kaka}) for $C\ne 0$.

\section{Appendix: Quasimeasures and multiplicators}
\label{sect-quasimeasures}

\subsection*{Quasimeasures}

Let us denote
by $\check g$
the inverse Fourier transform
of a tempered distribution $g$:
\[
\check g(t)=\mathcal{F}\sb{\omega\to t}^{-1}[g(\omega)].
\]

\begin{definition}\label{def-quasimeasure}

A tempered distribution
$\mu(\omega)$ is a
\emph{quasimeasure}
if  $\check \mu\in C\sb{b}(\R)$.
\end{definition}

For example,
any function from $L\sp 1(\R)$ is a quasimeasure,
and so is any finite Borel measure on $\R$.

\begin{lemma}\label{qmb}
Let $\mu(\omega)$ be a quasimeasure and $\varphi(\omega)$ be a test function
from the Schwartz space $\mathscr{S}(\R)$.
Then
\begin{equation}\label{mup}
\abs{\langle \mu(\omega),\varphi(\omega)\rangle}\le
C
\Vert
\check\varphi(t)
\Vert\sb{L^1(\R)}.
\end{equation}
\end{lemma}

The lemma is a trivial consequence of the Parseval identity:
\begin{equation}\label{mupp}
\abs{\langle \mu(\omega),\varphi(\omega)\rangle}
=
{2\pi}\abs{\langle \check\mu(t),\check\varphi(t)\rangle}
\le
{2\pi}
\Vert
\check\mu(t)
\Vert\sb{L^\infty(\R)}
\Vert
\check\varphi(t)
\Vert\sb{L^1(\R)}.
\end{equation}

\begin{definition}\label{dA}
$C\sb{b,F}(\R)$
is the vector space of bounded  functions
$f(t)\in C\sb{b}(\R)$
endowed with the following convergence:
$f\sb\varepsilon(t)\stackrel{C\sb{b,F}}\longrightarrow f(t)$,
$\varepsilon\to 0+\,$ if and only if
\begin{enumerate}
\item
$\forall T>0$,
\quad
$\Vert f\sb\varepsilon(t)-f(t)\Vert\sb{C[-T,T]}
\to 0$, $\varepsilon\to 0+$;
\item
$\sup\limits\sb{\varepsilon\in (0,1]}
\Vert f\sb\varepsilon(t)\Vert\sb{C\sb{b}(\R)}
<\infty$.

\end{enumerate}

\end{definition}

This type of convergence
coincides with the convergence
stated in the Ascoli-Arzel\`a theorem.
Next we introduce the dual class of the
``Ascoli-Arzel\`a quasimeasures''.
\begin{definition}\label{dQ}
$\mathcal{QM}(\R)$
is the
linear space of all quasimeasures
$\mu(\omega)$ endowed with the following convergence:
\[
\mu\sb\varepsilon(\omega)\stackrel{\mathcal{QM}}
{\mathop{\longrightarrow}
\limits\sb{\varepsilon\to 0+}}
\mu(\omega)
\quad{\rm if\ and\ only\ if}\quad
\check\mu\sb\varepsilon(t)\stackrel{C\sb{b,F}}
{\mathop{\longrightarrow}\limits\sb{\varepsilon\to 0+}}
\check\mu(t).
\]
\end{definition}

\subsection*{Multiplicators}

Now let us give a simple characterization
of multiplicators in $\mathcal{QM}(\R)$.
Let us consider a continuous function $M(\omega)\in C(\R)$.
We also denote by $M$ the
corresponding operator of multiplication:
\[
M:\;\mu(\omega)\mapsto M(\omega)\mu(\omega),
\qquad
\mu(\omega)\in C\sb{0}^\infty(\R).
\]

\begin{lemma}\label{lemma-quasi-1}
\begin{enumerate}
\item
Let
$\check M(t)\in L\sp 1(\R)$.
Then the operator $M$
extends
to a linear continuous operator in the space of quasimeasures:
\[
M:\;\mathcal{QM}(\R)\to\mathcal{QM}(\R).
\]
\item
Let $\mu\sb\varepsilon(\omega)\stackrel{\mathcal{QM}}\longrightarrow \mu(\omega)$  and
$\check M\sb\varepsilon(t)\stackrel{L\sp 1}\longrightarrow\check M(t)$
as $\varepsilon\to 0+$.
Then
\begin{equation}\label{mve}
M\sb\varepsilon(\omega)\mu\sb\varepsilon(\omega)
\stackrel{\mathcal{QM}}\longrightarrow M(\omega)\mu(\omega),
\qquad\varepsilon\to 0+.
\end{equation}
\end{enumerate}
\end{lemma}
\begin{proof}
First we define
$M(\omega)\mu(\omega):=\mathcal{F}\sb{t\to\omega}
[\big(\check M\ast\check\mu\big)(t)](\omega)$
that agrees with the case
$\mu\in C\sb{0}^\infty(\R)$.
Then
({\it i}) follows from ({\it ii})
with
$M\sb\varepsilon=M$
and $\mu\sb\varepsilon\in C\sb{0}^\infty(\R)$.
To prove  ({\it ii}),
we need to show that
\begin{equation}\label{lemma-second-statement}
\mathcal{F}\sb{\omega\to t}^{-1}[M\sb\varepsilon(\omega)\mu\sb\varepsilon(\omega)]
=\big(\check M\sb\varepsilon\ast\check\mu\sb\varepsilon\big)(t)
\stackrel{C\sb{b,F}}\longrightarrow
\big(\check M\ast\check\mu\big)(t).
\end{equation}
We have to check both conditions ({\it i}) and ({\it ii})
of Definition~\ref{dA} for the functions
\begin{eqnarray}
&&
f\sb\varepsilon(t):=\mathcal{F}\sb{\omega\to t}^{-1}
[M\sb\varepsilon(\omega)\mu\sb\varepsilon(\omega)]
=\big(\check M\sb\varepsilon\ast\check\mu\sb\varepsilon\big)(t),
\phantom{\int}
\nonumber
\\
&&
f(t):=\mathcal{F}\sb{\omega\to t}^{-1}
[M(\omega)\mu(\omega)]=
\big(\check M\ast\check\mu\big)(t).
\nonumber
\end{eqnarray}
We have:
\[
f\sb\varepsilon(t)-f(t)=
\big(\check M\sb\varepsilon\ast\check\mu\sb\varepsilon\big)(t)
-\big(\check M\ast\check\mu\big)(t)
=
\big(\!(\check M\sb\varepsilon-\check M)\ast\check\mu\sb\varepsilon\big)(t)
+\big(\check M\ast(\check\mu\sb\varepsilon-\check\mu)\!\big)(t).
\]
The first term in the right-hand side converges to zero
uniformly in $t\in\R$
since
 $\check M\sb\varepsilon-\check M\to 0$ in $L\sp 1$ while
$\check\mu\sb\varepsilon\in C\sb{b}(\R)$ and is bounded uniformly
for $\varepsilon\in(0,1)$.
Let us analyze the second term,
\begin{equation}\label{lim-second-term}
\int\sb{\R}\check M(\tau)(\check\mu\sb\varepsilon(t-\tau)-\check\mu(t-\tau))\,d\tau.
\end{equation}
Since $\check M\in L\sp 1$,
for any $\delta>0$
there exists a finite $R>0$ so that
$
\int\sb{\abs{\tau}>R}
\abs{\check M(\tau)}\,d\tau
\le\delta.
$
On the other hand, for any $T>0$,
the difference $\check\mu\sb\varepsilon(t-\tau)-\check\mu(t-\tau)$
is uniformly small for  $\abs{t}\le T$,
$\abs{\tau}<R$ and small $\varepsilon$.
Therefore,
the integral  (\ref{lim-second-term})
converges to zero uniformly in
$\abs{t}\le T$ as $\varepsilon\to 0+$.
Hence, the convergence  ({\it i}) of Definition~\ref{dA} follows.

Finally, the uniform bound ({\it ii}) of Definition~\ref{dA}
for the functions $f\sb\varepsilon(t)$ is obvious.
The convergence (\ref{lemma-second-statement}) is proved.
\qed\end{proof}

\subsection*{Bounds for multiplicators}

Let us justify the properties of the multiplicators
which we used in Section~\ref{sect-bound}.
Recall that we use the notation
\[
M\sb{x,\varepsilon}(\omega)
:=e\sp{ik(\omega+i\varepsilon)\abs{x}}\zeta(\omega),
\qquad
x\in\R,\quad \varepsilon\ge 0,
\]
where $\zeta(\omega)\in C\sb 0\sp\infty(\R)$
is a fixed cutoff function,
and also the notation
\[
N\sb x(\omega)
:=(ik(\omega)\sgn x)\sp\xx e\sp{ik(\omega)\abs{x}}
(-i\omega)\sp\yy \zeta(\omega),
\qquad
x\in\R,
\]
where
$\xx$, $\yy$ are fixed nonnegative integers.

\begin{lemma}\label{lemma-quasi-2}
For any fixed $x\in\R$ we have:
\begin{enumerate}
\item
$\check M\sb{x,\varepsilon}(t)\in L\sp 1(\R)$ for any
$\varepsilon\ge 0$.
\item
$\check M\sb{x,\varepsilon}(t)\stackrel{L\sp 1}
\longrightarrow\check M\sb{x,0}(t)$,
\quad$\varepsilon\to 0$.
\item
$\check N\sb x \in L\sp 1(\R)$, and for any $R>0$
there exists $C\sb{\xx,\yy,R}>0$
so that
\begin{equation}\label{B*}
\sup\limits\sb{\abs{x}\le R}
\,
\norm{
\check N\sb x}\sb{L\sp 1(\R)}\le C\sb{\xx,\yy,R}.
\end{equation}
\end{enumerate}
\end{lemma}

\begin{proof}
For any fixed $x\in\R$, the Puiseux expansion holds:
\begin{equation}\label{MHfA}
e\sp{ik(\omega+i\varepsilon)\abs{x}}
\sim 1+
\sum\sb{\pm}\sum\limits\sb{j=1}^{\infty}
C\sb j\sp\pm
(x)(\omega+i\varepsilon\mp m)\sp{j/2},
\qquad
\omega+i\varepsilon\to \pm m,
\qquad
\varepsilon >0.
\end{equation}
Therefore,
the function $\check M\sb{x,\varepsilon}(t)$
is smooth and
decays at least like $\abs{t}\sp{-3/2}$
when $t\to\infty$.
This finishes the proof of the first statement of the lemma.

The second statement of the lemma follows from (\ref{MHfA}).

The last statement of the lemma
follows by the same arguments from
the Puiseux expansion for $\check N\sb x(\omega)$
similar to expansion (\ref{MHfA}) with $\varepsilon=0$.
\qed\end{proof}

\section{Appendix: Global well-posedness}
\label{sect-existence}

Here we prove Theorem~\ref{theorem-well-posedness}.
We first need
to adjust the nonlinearity $F$
so that it becomes bounded, together with
its derivatives.
Define
\begin{equation}\label{def-Lambda}
\Lambda(\Psi\sb 0)
=\sqrt{\frac{\mathcal{H}(\Psi\sb 0)-{A}}{m-{B}}},
\end{equation}
where $\Psi\sb 0\in{\E}$ is the initial data
from Theorem~\ref{theorem-well-posedness}
and ${A}$, ${B}$ are constants from (\ref{bound-below}).
Then we may pick a modified potential function
$\widetilde{U}\in C\sp 2(\C,\R)$,
$\widetilde{U}(\psi)=\widetilde{U}(\abs{\psi})$,
so that
\begin{equation}\label{new-U}
\widetilde{U}(\psi)=U(\psi)
\qquad{\rm for}\ \abs{\psi}\le \Lambda(\Psi\sb 0),
\quad
\psi\in\C,
\end{equation}
$\widetilde{U}(\psi)$ satisfies (\ref{bound-below})
with the same constants ${A}$, ${B}$ as $U(\psi)$ does:
\begin{equation}\label{new-U-2}
\widetilde{U}(\psi)\ge {A}-{B}\abs{\psi}^2,
\quad{\rm for}\ \psi\in\C,\quad
{\rm where}\ {A}\in\R\ {\rm and}\ 0\le {B}<m,
\end{equation}
and so that
$\abs{\widetilde{U}(\psi)}$,
$\abs{\widetilde{U}'(\psi)}$,
and $\abs{\widetilde{U}''(\psi)}$
are bounded for $\psi\ge 0$.
We define
\begin{equation}\label{f-reg}
\widetilde{F}(\psi)
=-\nabla \widetilde{U}(\psi),
\qquad\psi\in\C,
\end{equation}
where $\nabla$ denotes the gradient with respect to
$\Re\psi$, $\Im\psi$;
Then $\widetilde{F}(e\sp{is}\psi)=e\sp{is}\widetilde{F}(\psi)$ for any $\psi\in\C$, $s\in\R$.

We consider the Cauchy problem of type
(\ref{KG}) with the modified
nonlinearity,
\begin{equation}\label{KG-a}
\left\{
\begin{array}{l}
\ddot\psi(x,t)
=\psi''(x,t)-m^2\psi(x,t)+\delta(x)\widetilde F(\psi(0,t)),
\qquad
x\in\R,
\quad t\in\R,
\\
\psi\at{t=0}=\psi\sb 0(x),
\qquad
\dot\psi\at{t=0}=\pi\sb 0(x),
\end{array}\right.
\end{equation}
which we rewrite in the vector form
in terms of
$\Psi=\left[\begin{array}{c}\psi(x,t)\\\pi(x,t)\end{array}\right]$,
similarly to (\ref{KG-cp}):
\begin{equation}
\dot\Psi
=\left[\begin{array}{cc}0&1\\\p\sb x^2-m^2&0\end{array}\right]\Psi
+\delta(x)\left[\begin{array}{c}0\\\widetilde{F}(\psi)\end{array}\right],
\quad
\quad\Psi\at{t=0}
=\Psi\sb 0
\equiv\left[\begin{array}{c}\psi\sb 0(x)\\\pi\sb 0(x)\end{array}\right].
\label{KG-ap}
\end{equation}
This is a Hamiltonian system, with the Hamilton functional
\begin{equation}\label{KG-a-h}
\widetilde{\mathcal{H}}(\Psi)
=\int\limits\sb\R
\left(
\abs{\pi}\sp 2
+\abs{\nabla\psi}\sp 2+m^2\abs{\psi}\sp 2
\right)\,dx
+\widetilde{U}(\psi(0,t)),
\quad
\Psi=\left[\begin{array}{c}\psi(x)\\
\pi(x)\end{array}\right]\in\E,
\end{equation}
which is Fr\'echet differentiable in the space
${\E}=H\sp 1\times L\sp 2$.
By the Sobolev embedding theorem,
$L\sp\infty(\R)\subset H\sp 1(\R)$,
and there is the following inequality:
\begin{equation}\label{sobolev-embedding}
\norm{\psi}\sb{L\sp\infty}^2
\le\frac{1}{2m}
(\norm{\psi'}\sb{L\sp 2}^2+m^2\norm{\psi}\sb{L\sp 2}^2)
\le\frac{1}{2m}
\norm{\Psi}\sb{\E}^2.
\end{equation}
Thus, (\ref{new-U-2}) leads to
\begin{equation}\label{bound-on-u}
\widetilde{U}(\psi(0))
\ge {A}-{B}\norm{\psi}\sb{L\sp\infty}^2
\ge {A}-\frac{{B}}{2m}\norm{\Psi}\sb{\E}^2.
\end{equation}
Taking into account (\ref{KG-a-h}),
we obtain the inequality
\begin{equation}
\norm{\Psi}\sb{\E}\sp 2
=2\widetilde{\mathcal{H}}(\Psi)-2 \widetilde{U}(\psi(0))
\le 2\widetilde{\mathcal{H}}(\Psi)-2{A}+\frac{{B}}{m}
\norm{\Psi}\sb{\E}^2,
\qquad
\Psi\in\E,
\end{equation}
which implies
\begin{equation}\label{t-bound-1}
\norm{\Psi}\sb{\E}\sp 2
\le\frac{2m}{m-{B}}\left(\widetilde{\mathcal{H}}(\Psi)-{A}\right),
\qquad\Psi\in\E.
\end{equation}

\begin{lemma}\label{lemma-same-u}
\begin{enumerate}
\item
There is the identity
$\widetilde{\mathcal{H}}(\Psi\sb 0)
=\mathcal{H}(\Psi\sb 0)$.
\item
If $\Psi=\left[\begin{array}{c}\psi(x)\\\pi(x)\end{array}\right]\in{\E}$
satisfies
$\widetilde{\mathcal{H}}(\Psi)\le\widetilde{\mathcal{H}}(\Psi\sb 0)$,
then
$\ \widetilde{U}(\psi(0))=U(\psi(0))$.
\end{enumerate}
\end{lemma}

\begin{proof}
\begin{enumerate}
\item
According to (\ref{t-bound-1}),
the Sobolev embedding (\ref{sobolev-embedding}),
and the choice of $\Lambda(\Psi\sb 0)$ in (\ref{def-Lambda}),
\begin{equation}
\norm{\psi\sb 0}\sb{L\sp\infty}^2
\le\frac{1}{2m}\norm{\Psi\sb 0}\sb{\E}^2
\le\frac{\mathcal{H}(\Psi\sb 0)-{A}}{m-{B}}
=\Lambda(\Psi\sb 0)^2.
\end{equation}
Thus, according to the choice of $\widetilde{U}$
(equality (\ref{new-U})),
$
\widetilde{U}(\psi\sb 0(0))=U(\psi\sb 0(0)),
$
proving ({\it i}).
\item
By (\ref{sobolev-embedding}), (\ref{t-bound-1}),
the condition
$\widetilde{\mathcal{H}}(\Psi)\le\widetilde{\mathcal{H}}(\Psi\sb 0)$,
and part ({\it i}) of the Lemma, we have:
\[
\norm{\psi}\sb{L\sp\infty}^2
\le\frac{1}{2m}\norm{\Psi}\sb{\E}^2
\le\frac{\widetilde{\mathcal{H}}(\Psi)-{A}}{m-{B}}
\le\frac{\widetilde{\mathcal{H}}(\Psi\sb 0)-{A}}{m-{B}}
=\frac{\mathcal{H}(\Psi\sb 0)-{A}}{m-{B}}
=\Lambda(\Psi\sb 0)^2.
\]
Hence, ({\it ii}) follows by (\ref{new-U}).
\end{enumerate}
\qed\end{proof}

\begin{remark}
We will show that if $\Psi(t)$ solves (\ref{KG-ap}),
then
$\widetilde{\mathcal{H}}(\Psi(t))=\widetilde{\mathcal{H}}(\Psi\sb 0)$,
and therefore
$\widetilde{U}(\psi(0,t))=U(\psi(0,t))$
by Lemma~\ref{lemma-same-u}~({\it ii}).
Hence, $\widetilde{F}(\psi(0,t))=F(\psi(0,t))$
for all $t\ge 0$,
allowing us to conclude that
$\Psi(t)$ solves (\ref{KG-cp})
as well as (\ref{KG-ap}).
\end{remark}

\subsection*{Local well-posedness}
The solution to the Cauchy problem
\begin{equation}
\dot\Xi
=\left[\begin{array}{cc}0&1\\
\p\sb x^2-m^2&0\end{array}\right]\Xi,
\qquad
\Xi(x,0)
=\Xi\sb 0(x)
=\left[\begin{array}{c}\xi\sb 0(x)\\\eta\sb 0(x)\end{array}\right]
\end{equation}
is represented by
\begin{equation}\label{integral-representation-free}
\Xi(x,t)
=W\sb 0(t)
\Xi\sb 0
=
\int\sb\R
\left[\begin{array}{cc}\dot{G}(x-y,t)&{G}(x-y,t)\\
\ddot{G}(x-y,t)&\dot{G}(x-y,t)\end{array}\right]
\left[\begin{array}{c}\xi\sb 0(y)\\\eta\sb 0(y)\end{array}\right]
\,dy,
\end{equation}
where ${G}(x,t)$
is the forward fundamental solution to the Klein-Gordon equation,
${G}(x,t)
=\theta(t-\abs{x})J\sb 0(m\sqrt{t\sp 2-x\sp 2})/2,
$
with $J\sb 0$ being the Bessel function
(see e.g. \cite{MR1364201}).
Then
the solution to the Cauchy problem (\ref{KG-ap})
can be represented by
\begin{eqnarray}\label{integral-representation}
&&
\Psi(x,t)
=W\sb 0(t)\Psi\sb 0
+{Z}[\psi(0,\cdot)](t),
\nonumber
\\
&&
{\rm where}
\quad
{Z}[\psi(0,\cdot)](t):=
\int\sb 0\sp{t}
W\sb 0(t-s)
\left[\begin{array}{c}0
\\
\delta(\cdot)\widetilde{F}(\psi(0,s))
\end{array}\right]
\,ds.
\end{eqnarray}

\begin{lemma}\label{lemma-bessel}
For any
nonnegative integers $\xx$ and $\yy$
there is a constant $C\sb{\xx,\yy}>0$ such that
\begin{equation}\label{bessel-many}
\abs{\p\sb x\sp\xx \p\sb t\sp\yy J\sb 0(m\sqrt{t^2-x^2})}
\le C\sb{\xx,\yy}(1+t)^{\xx+\yy},
\qquad
\abs{x}<t.
\end{equation}
\end{lemma}

\begin{proof}
The proof immediately follows from the observation
that all the derivatives of the Bessel function  $J\sb 0(z)$
are bounded for $z\in\R$, and that $J\sb 0(z)$
is an absolutely converging Taylor series in even powers of $z$.
Hence,
all derivatives of the function $J\sb 0(\sqrt{r})$
in $r$ are bounded for $r\ge 0$.
\qed\end{proof}

The next lemma establishes
the contraction principle
for the integral equation
(\ref{integral-representation}).

\begin{lemma}\label{lemma-bounds}
There exists a constant $C>0$ so that
for any two functions
$
\Psi\sb{k}(\cdot,t)
=
\left[\begin{array}{c}\psi\sb{k}(\cdot,t)\\
\pi\sb{k}(\cdot,t)
\end{array}\right]\in C([0,1],\E),
$
$k=1,\,2$,
we have:
\[
\norm{
{Z}[\psi\sb 1(0,\cdot)](t)-{Z}[\psi\sb 2(0,\cdot)](t)
}\sb{\E}
\le
C t^{1/2}
\sup\sb{s\in[0,t]}\norm{\Psi\sb 1(\cdot,s)-\Psi\sb 2(\cdot,s)}\sb{\E},
\]
for $0\le t\le 1$.
\end{lemma}

\begin{proof}
According to (\ref{integral-representation-free}) and
(\ref{integral-representation}),
\[
{Z}[\psi\sb 1(0,\cdot)](t)-{Z}[\psi\sb 2(0,\cdot)](t)=
\left[\begin{array}{c}
I(x,t)\\
\partial_t I(x,t)
\end{array}\right],
\]
where
\[
I(x,t):=\displaystyle\int\sb 0\sp t
{G}(x,t-s)
\left(
\widetilde{F}(\psi\sb 1(0,s))-\widetilde{F}(\psi\sb 2(0,s))
\right)\,ds.
\]
First we prove the
$L\sp 2$ estimate for $I(x,t)$.
By the Sobolev embedding theorem,
\begin{eqnarray}
\Norm{
I(\cdot,t)
}\sb{L\sp 2}
&\le&
\cnst
\Norm{
\int\sb 0\sp t
\theta(t-s-\abs{x})
\abs{\widetilde{F}(\psi\sb 1(0,s))-\widetilde{F}(\psi\sb 2(0,s))}
\,ds
}\sb{L\sp 2}
\nonumber\\
\nonumber\\
&
\le&
\cnst
\sup\sb{z\in\C}\abs{\nabla\widetilde{F}(z)}
\Norm{
\int\sb 0\sp t
\theta(t-s-\abs{x})
\,ds
}\sb{L\sp 2}
\,
\sup\sb{s\in[0,t]}
\norm{\psi\sb 1(\cdot,s)-\psi\sb 2(\cdot,s)}\sb{H\sp 1}
\nonumber\\
\nonumber\\
&
\le&
C'
\,t\sp{3/2}
\sup\sb{s\in[0,t]}
\norm{\psi\sb 1(\cdot,s)-\psi\sb 2(\cdot,s)}\sb{H\sp 1},
\label{e-l2}
\end{eqnarray}
where we took into account that
$\abs{\nabla\widetilde{F}(z)}$
is bounded
due to the choice of $\widetilde{U}$.

Similarly, we derive
the  $L\sp 2$ estimate
for the derivative $\p\sb x I(x,t)$.
We first analyze
\[
\p\sb x{G}(x,t-s)=
\frac 1 2
\theta(t-s-\abs{x})\p\sb x J\sb 0(m\sqrt{(t-s)\sp 2-x\sp 2})
-\frac 1 2\sgn{x}\,\delta(t-s-\abs{x}).
\]
By Lemma~\ref{lemma-bessel} for $\abs{x}\le\abs{t-s}\le 1$,
we have
$
\abs{
\p\sb x J\sb 0(m\sqrt{(t-s)\sp 2-x\sp 2})}
\le\cnst;
$
We conclude that
$
\Norm{
\p\sb x I(\cdot,t)
}\sb{L^2}
$
is bounded by
\begin{eqnarray}\label{d-e-l2}
&&
\Norm{
\int\sb 0\sp t
\left[C
\theta(t-s-\abs{x})
+
\frac{\delta(t-s-\abs{x})}{2}
\right]
\,ds
}\sb{L\sp 2}
\,\sup\sb{s\in[0,t]}
\Abs{\widetilde{F}(\psi\sb 1(0,s))-\widetilde{F}(\psi\sb 2(0,s))}
\nonumber\\
\nonumber\\
&&
\qquad\qquad
\le
C
\Norm{
\theta(t-\abs{x})
}\sb{L\sp 2}
\,\sup\sb{s\in[0,t]}
\norm{\psi\sb 1(\cdot,s)-\psi\sb 2(\cdot,s)}\sb{H\sp 1}
\nonumber\\
\nonumber\\
&&
\qquad\qquad
\le
C'
\,t^{1/2}
\,\sup\sb{s\in[0,t]}
\norm{\psi\sb 1(\cdot,s)-\psi\sb 2(\cdot,s)}\sb{H\sp 1}.
\end{eqnarray}
The $L^2$-norm of $\p\sb t I(x,t)$
is estimated similarly.
\qed\end{proof}

For $E>0$,
let us denote
${\E}\sb E=\{\Psi\sb 0\in \E
\sothat
\mathcal{H}(\Psi\sb 0)\le E\}$.

\begin{corollary}
\label{cor-l8-existence}
\begin{enumerate}
\item
For any $E>0$ there exists $\tau=\tau(E)>0$
such that for any $\Psi\sb 0\in {\E}\sb E$
there is a unique solution
$\Psi(x,t)\in C([0,\tau],{\E})$
to the Cauchy problem {\rm (\ref{KG-ap})}
with the initial condition
 $\Psi(0)=\Psi\sb 0$.
\item
The map $W(t):\;\Psi\sb 0\mapsto \Psi(t)$,
 $t\in[0,\tau]$
are continuous maps from ${\E}\sb E$ to $\E$.
\end{enumerate}
\end{corollary}

\subsection*{Smoothness of the solution}

In this section,
we will study the smoothness of the solution
\[
\Psi(x,t)=(\psi(x,t),\pi(x,t))\in C([0,\tau],\E)
\]
constructed in Corollary~\ref{cor-l8-existence}~({\it i})
assuming that $\psi\sb 0(x),\pi\sb 0(x)\in C\sb 0^\infty(\R)$.
According to the integral representation (\ref{integral-representation}),
$\psi(x,t)$, $t\in [0,\tau]$, can be represented as
\begin{equation}\label{integral-representation-psi}
\psi(x,t)
=
\int\limits\sb\R
\left(
\dot{G}(x-y,t)\psi\sb 0(y)+{G}(x-y,t)\pi\sb 0(y)
\right)
dy
+
\int\limits\sb 0\sp t {G}(x,t-s)\widetilde{F}(\psi(0,s))\,ds.
\end{equation}

First, let us prove the smoothness
of the function $\psi(0,t)$.

\begin{lemma}\label{lemma-smoothness-0}
$\psi(0,t)\in C\sp\infty([0,\tau])$.
\end{lemma}

\begin{proof}
The integral representation (\ref{integral-representation-psi})
implies  that,
for $t\in [0,\tau]$,
\begin{equation}\label{imt}
\psi(0,t)
=
\int\limits\sb\R
\left(
\dot{G}(y,t)\psi\sb 0(y)+{G}(y,t)\pi\sb 0(y)
\right)dy
+
\frac{1}{2}
\int\limits\sb 0\sp t
J\sb 0(m(t-s))
\widetilde{F}(\psi(0,s))
\,ds.
\end{equation}
The first integral is a smooth function.
Further,
from $\norm{\psi(\cdot,t)}\sb{H\sp 1}\le C<\infty$, $t\in [0,\tau]$,
we conclude that $\abs{\psi(0,t)}$ is bounded.
Hence, (\ref{imt}) implies
that $\psi(0,\cdot)\in C([0,\tau])$,
and then by induction that $\psi(0,\cdot)\in C\sp\infty([0,\tau])$
since the Bessel function is smooth.
\qed\end{proof}

Now, from (\ref{integral-representation-psi}),
we conclude that $\psi(x,t)$ is smooth
away from the singularities of ${G}(x,t)$.

\begin{proposition}\label{coro}
The solution
$\psi(x,t)$ is piecewise smooth
inside each of the four regions of $[0,\tau]\times\R$
cut off by the lines $x=0$ and $x=\pm t$.
\end{proposition}

\begin{proof}
The first integral
in the right-hand side of (\ref{integral-representation-psi})
is infinitely smooth in $x$ and $t$ for all $x\in\R$, $t\ge 0$.
Now let us consider the second integral
in the right-hand side of (\ref{integral-representation-psi}),
which could be written as follows:
\begin{equation}\label{ir-psi-two}
\frac{\theta(t-\abs{x})}{2}
\int\sb 0\sp{t-\abs{x}}
J\sb 0(m\sqrt{(t-s)^2-x^2})
\widetilde{F}(\psi(0,s))\,ds.
\end{equation}
Here the function
$\widetilde{F}(\psi(0,s))$
is smooth in $s\in[0,\tau]$ by Lemma~\ref{lemma-smoothness-0}.
By Lemma~\ref{lemma-bessel},
all the partial derivatives of
$J\sb 0(m\sqrt{(t-s)^2-x^2})$
in $x$ and $t$
are continuous and uniformly bounded for $\abs{x}<t-s$,
$t\le\tau$.
Therefore, (\ref{ir-psi-two})
is smooth,
with all the derivatives uniformly bounded,
in each of the regions
$0\le x\le t$, $-t\le x\le 0$.
In the regions $\abs{x}>t$,
(\ref{ir-psi-two}) is identically equal to zero.
\qed\end{proof}

\begin{lemma}\label{lemma-psi-dot}
For
$0<t\le\tau$,
\begin{equation}\label{psi-dot-pm}
\lim\sb{x\to 0-}\dot\psi(x,t)=\lim\sb{x\to 0+}\dot\psi(x,t).
\end{equation}
\end{lemma}

\begin{proof}
We have to analyze only
the contribution from the second term
in the right-hand side of (\ref{integral-representation-psi}),
that is,
\[
\p\sb t
\int\limits\sb 0\sp t {G}(x,t-s)\widetilde{F}(\psi(0,s))\,ds
=
{G}(x,0+)\widetilde{F}(\psi(0,t))
+
\int\limits\sb 0\sp t
\dot{G}(x,t-s)\widetilde{F}(\psi(0,s))\,ds.
\]
The first term in the right-hand side
is equal to zero for $x\ne 0$.
The second term
is continuous since
the Green function
${G}(x,t-s)$
is smooth at $x=0$ for $t-s>0$.
\qed\end{proof}

\begin{lemma}\label{lemma-psi-dot-prime}
For
$0<t\le\tau$,
\begin{enumerate}
\item
$
\dot\psi(x,t)+\psi'(x,t)
$
is continuous across
the characteristic $x=t$.
\item
$
\dot\psi(x,t)-\psi'(x,t)
$
is continuous across the characteristic $x=-t$.
\end{enumerate}
\end{lemma}

\begin{proof}
The proofs for both statements of the Lemma are identical;
we will only prove the first statement with $x>0$.
We have to study only
the contribution from the second term
in the right-hand side of (\ref{integral-representation-psi}),
i.e.
\begin{equation}\label{psi-dot-rhs-2}
(\p\sb t+\p\sb x)
\int\limits\sb 0\sp t {G}(x,t-s)\widetilde{F}(\psi(0,s))\,ds
=
\int\limits\sb 0\sp t
(\p\sb t+\p\sb x){G}(x,t-s)
\widetilde{F}(\psi(0,s))\,ds.
\end{equation}
Here we took into account that, as above,
${G}(x,0+)\widetilde{F}(\psi(0,t))=0$ for $x\ne 0$.
Next key observation
is that,
for $x>0$, the derivative $\p\sb t+\p\sb x$
applied to ${G}(x,t)$,
does not produce a delta-function:
\[
(\p\sb t+\p\sb x){G}(x,t)
=\frac 1 2
\left\{
\theta(t-x)
(\p\sb t+\p\sb x)J\sb 0(m\sqrt{t^2-x^2})
\right\}.
\]
Hence, the integral (\ref{psi-dot-rhs-2})
is continuous in $x$ and $t$ across the line $x=t$, $0<t\le\tau$
by Lemma~\ref{lemma-bessel}
\qed\end{proof}

\subsection*{Energy conservation and global well-posedness}

\begin{lemma}\label{lemma-e-conserved}
For the solution
to the Cauchy problem (\ref{KG-ap})
with the initial data
$\Psi\sb 0\in\E$,
the energy is conserved:
$\widetilde{\mathcal{H}}(\Psi(t))=\const$, $t\in [0,\tau]$.
\end{lemma}

\begin{proof}
We follow \cite{MR1359949}.
First, we prove that the energy is conserved
for the smooth initial data with compact support:
$\Psi\sb 0=\left[\begin{array}{c}\psi\sb 0\\\pi\sb 0\end{array}\right]$,
with
$\psi\sb 0$, $\pi\sb 0\in C\sp\infty\sb 0(\R)$.
Consider the norm (\ref{def-e}),
\begin{equation}
\norm{\Psi(t)}\sb{\E}^2
=
\int\sb{-\infty}\sp\infty
[
\abs{\dot\psi}\sp 2+\abs{\psi'}\sp 2+m^2\abs{\psi}\sp 2
]\,dx,
\qquad
t\in[0,\tau].
\end{equation}
We split this integral into four pieces: The integration over
$(-\infty,-t)$, $(-t,0)$, $(0,t)$, and $(t,\infty)$.
By Proposition~\ref{coro},
on the support of each of these integrals
$\psi(x,t)$ for $t\in[0,\tau]$
is a smooth function of $x$ and $t$.
Then, differentiating,
we may express $\p\sb t \norm{\Psi(t)}\sb{\E}^2$ as
\begin{eqnarray}\label{p-t-e}
&&
\p\sb t \norm{\Psi(t)}\sb{\E}^2
=
\Big[
\abs{\dot\psi}\sp 2+\abs{\psi'}\sp 2+m^2\abs{\psi}\sp 2
\Big]\sb{x=-t-0}\sp{x=-t+0}
\\
&&
-
\Big[
\abs{\dot\psi}\sp 2+\abs{\psi'}\sp 2+m^2\abs{\psi}\sp 2
\Big]\sb{x=t-0}\sp{x=t+0}
+2\int\sb{-\infty}\sp{\infty}
[
\dot\psi\ddot\psi+\psi'\dot\psi'+m^2\psi\dot\psi
]\,dx,
\quad
t\in[0,\tau].
\nonumber
\end{eqnarray}
The terms $m^2\abs{\psi}\sp 2$ could be discarded
due to continuity of $\psi$
across the characteristics $x=\pm t$.
Integrating by parts the terms $\psi'\dot\psi'$
and using the cancellations of the integrals
due to equation (\ref{KG-a}) away from $x=0$,
we get:
\begin{eqnarray}
&&\p\sb t \norm{\Psi(t)}\sb{\E}^2
=
\Big[
\abs{\dot\psi}\sp 2+\abs{\psi'}\sp 2
-2\psi'\dot\psi
\Big]\sb{x=-t-0}\sp{x=-t+0}
\nonumber\\
&&\qquad\qquad
-
\Big[
\abs{\dot\psi}\sp 2+\abs{\psi'}\sp 2
+2\psi'\dot\psi
\Big]\sb{x=t-0}\sp{x=t+0}
-2\Big[\psi'\dot\psi\Big]\sb{x=0-}\sp{x=0+}
\nonumber\\
\nonumber\\
&&\qquad
=
\Big[
(\dot\psi-\psi')\sp 2
\Big]\sb{x=-t-0}\sp{x=-t+0}
-\Big[
(\dot\psi+\psi')\sp 2
\Big]\sb{x=t-0}\sp{x=t+0}
-2\Big[\psi'\dot\psi\Big]\sb{x=0-}\sp{x=0+}.
\label{last-line}
\end{eqnarray}
According to Lemma~\ref{lemma-psi-dot-prime},
the first two terms in (\ref{last-line})
 do not give any contribution.
Let us compute the contribution of the last term.
According to Lemma~\ref{lemma-psi-dot},
$\dot\psi(0\pm,t)=\dot\psi(0,t)$
for $t\in [0,\tau]$, therefore
\[
\left[\psi'\dot\psi\right]\sb{x=0-}\sp{x=0+}
=\left[\psi'(x,t)\right]\sb{x=0-}\sp{x=0+}\psi(0,t)
=
-\widetilde{F}(\psi(0,t))\dot\psi(0,t)
=\frac{d}{dt}\widetilde{U}(\psi(0,t)).
\]
In the second equality,
we computed the jump of $\psi'$
using equation (\ref{KG-a}) and the piecewise smoothness of the solution.
We conclude that
\[
\frac{d}{dt}
\left\{\frac 1 2\norm{\Psi(t)}\sb{\E}^2+\widetilde{U}(\psi(0,t))\right\}=0,
\]
and hence the value of the functional
$\widetilde{\mathcal{H}}$
defined in (\ref{KG-a-h}) is conserved.

Since we proved the energy conservation
for the initial data that constitute a dense set
in ${\E}$
and since the dynamical group
is continuous
in ${\E}$
by
Corollary~\ref{cor-l8-existence}~({\it ii}),
we conclude that the energy is conserved
for arbitrary initial data from ${\E}$.
\qed\end{proof}

\begin{corollary}
\label{coroll}
\begin{enumerate}
\item
The solution $\Psi$ to the Cauchy problem (\ref{KG-ap})
with the initial data
$\Psi\at{t=0}=\Psi\sb 0\in {\E}$
exists globally:
$
\Psi\in C\sb{b}(\R,{\E}).
$
\item
The energy
is conserved:
$
\widetilde{\mathcal{H}}(\Psi(t))=\widetilde{\mathcal{H}}(\Psi\sb 0),
\qquad t\ge 0.
$
\end{enumerate}
\end{corollary}

\begin{proof}
Corollary~\ref{cor-l8-existence}~({\it i})
yields a solution
$\Psi\in L\sp\infty([0,\tau],{\E})$
with a positive $\tau=\tau(E)$.
However,
the value of $\mathcal{H}(\Psi(t))$
is conserved for $t\le \tau$
by Lemma~\ref{lemma-e-conserved}.
Corollary~\ref{cor-l8-existence}~({\it i})
allows us to extend $\Psi$
to the interval $[\tau,2 \tau]$,
and eventually to all $t\ge 0$.
In the same way we extend the solution $\Psi(t)$
for all $t<0$.
\qed\end{proof}

\subsection*{Conclusion of the proof of Theorem~\ref{theorem-well-posedness}}

The trajectory
$\Psi
=\scriptsize{\left[\!\begin{array}{c}\psi(x,t)\\\pi(x,t)\end{array}\!\right]}
\in C\sb{b}(\R,{\E})$
is a solution to (\ref{KG-ap}),
for which
Corollary~\ref{coroll}~({\it ii})
together with Lemma~\ref{lemma-same-u}~({\it i})
imply the energy conservation  (\ref{ec}).
By Lemma~\ref{lemma-same-u}~({\it ii}),
$\widetilde{U}(\psi(0,t))=U(\psi(0,t))$,
for all $t\in\R$.
This tells us that
$\psi(x,t)$ is a solution to (\ref{KG}).
Finally, the a priori bound (\ref{eb}) follows from
(\ref{t-bound-1}) and the conservation of $\mathcal{H}(\Psi(t))$.
This finishes the proof of Theorem~\ref{theorem-well-posedness}.

\noindent
{\bf Acknowledgements.}
The authors thank H.~Brezis,
V.S.~Buslaev, J.~Ginibre, P.-L. Lions, L.~Nirenberg, J.~Shatah,
A.~Shnirelman, H.~Spohn,
W.~Strauss, G.~Velo, and M.I.~Vishik for fruitful discussions,
and
A.~Merzon
for helpful remarks on the Titchmarsh theorem.
The authors are also indebted to P.~Joly, F.~Collino, and T.~Fouquet
from Project ONDES (INRIA),
and to A.~Vinnichenko
for the help with numerical experiments
for nonlinear wave equations.

%

\bibliographystyle{amsalpha}
\bibliography{ubk-mathsci,ubk-local}

\end{document}